\documentclass[11pt,leqno,twoside]{article}
\usepackage[utf8]{inputenc}
\usepackage[T1]{fontenc}
\usepackage[english]{babel}
\usepackage{authblk}
\usepackage{amsfonts}
\usepackage{tcolorbox}
\usepackage{amsmath}
\usepackage{amsthm}
\usepackage{geometry}
\usepackage{stmaryrd}
\usepackage{amssymb}
\usepackage{enumitem}
\usepackage{mathrsfs}
\usepackage{listings}
\usepackage{empheq,mathtools}
\usepackage{graphicx}
\usepackage{fullpage}
\usepackage{color}
\usepackage{hyperref}
\usepackage{fancyhdr}
\usepackage{bbm}
\usepackage{esint}
\usepackage{physics}
\usepackage{framed}
\usepackage{emptypage}
\usepackage{csquotes}
\usepackage[normalem]{ulem}
\usepackage[
backend=biber,
bibstyle=alphabetic,
citestyle=alphabetic,
giveninits=true,
sorting=nyt
]{biblatex}

\iffalse
\DeclareFieldFormat*{title}{\emph{#1\isdot}}					

\DeclareFieldFormat*{titlecase}{
    \ifdef{\currentfield}
      {\ifcurrentfield{title}
         {\usefield{\uline}{\currentfield}}%
         {#1}}
      {#1}}
\fi

\hypersetup{
	colorlinks=true,
	linkcolor=blue,
	filecolor=magenta,
	urlcolor=cyan,
	citecolor=blue,
	pdfpagemode=FullScreen,
}

\definecolor{darkWhite}{rgb}{0.94,0.94,0.94}
\definecolor{antiquebrass}{rgb}{0.8, 0.58, 0.46}
\definecolor{Bleu}{RGB}{23,132,255}

\definecolor{Defi}{RGB}{0,255,255}
\definecolor{Prop}{RGB}{0,255,0}
\definecolor{Lemm}{RGB}{255,255,0}
\definecolor{Theo}{RGB}{255,0,0}
\definecolor{Coro}{RGB}{255,0,255}

\definecolor{shadecolor}{RGB}{230, 230, 255}

\newtheorem*{remark}{Remark}
\newtheorem{thm}{Theorem}[section]
\newtheorem{prop}[thm]{Proposition}
\newtheorem{cor}[thm]{Corollary}
\newtheorem{lem}[thm]{Lemma}
\newtheorem{defprop}[thm]{Definition / Proposition}
\newtheorem{law}[thm]{Definition}

\newenvironment{Theoreme}{
	\begin{center}\begin{minipage}{\linewidth}\begin{thm}
}
{
	\end{thm}\end{minipage}\end{center}
}

\newenvironment{Proposition}{
	\begin{center}\begin{minipage}{\linewidth}
	\begin{prop}
}
{
	\end{prop}
	\end{minipage}\end{center}
}

\newenvironment{Corollaire}{
	\begin{center}\begin{minipage}{\linewidth}
	\begin{cor}
}
{
	\end{cor}
	\end{minipage}\end{center}
}

\newenvironment{Lemme}{
	\begin{center}\begin{minipage}{\linewidth}
	\begin{lem}
}
{
	\end{lem}
	\end{minipage}\end{center}
}

\newenvironment{Preuve}[1][\unskip]{
	\paragraph{Proof #1:}
}
{	
	\begin{flushright}
		\qedsymbol\smallskip
	\end{flushright}
}

\DeclarePairedDelimiterX\set[1]\lbrace\rbrace{\def\given{\;\delimsize\vert\;}#1}

\newcommand{\new}{\hfill \break \indent}

\newcommand{\npa}[1]{ {(#1)} }
\newcommand{\bpa}[1]{ {\big( #1 \big)} }

\newcommand{\lrpa}[1]{ {\mathopen{}\left( #1 \right)\mathclose{}} }

\newcommand{\bbr}[1]{ {\big[ #1 \big]} }

\newcommand{\Bbbr}[1]{ {\Big[ #1 \Big]} }

\newcommand{\lrbr}[1]{ {\left[ #1 \right]} }

\newcommand{\ga}{ {\alpha} }

\renewcommand{\gg}{ {\gamma} }
\newcommand{\gG}{ {\Gamma} }
\newcommand{\gd}{ {\delta} }
\newcommand{\gD}{ {\Delta} }

\newcommand{\gve}{ {\varepsilon} }

\newcommand{\gt}{ {\theta} }

\newcommand{\gk}{ {\kappa} }
\newcommand{\gl}{ {\lambda} }
\newcommand{\gL}{ {\Lambda} }

\newcommand{\gs}{ {\sigma} }
\newcommand{\gS}{ {\Sigma} }

\newcommand{\gO}{ {\Omega} }
\newcommand{\gvp}{ {\varphi} }

\newcommand{\NN}{ {\mathbb{N}} }

\newcommand{\RR}{ {\mathbb{R}} }

\renewcommand{\SS}{ {\mathbb{S}} }

\newcommand{\cA}{ {\mathcal{A}} }

\newcommand{\cC}{ {\mathcal{C}} }

\newcommand{\cF}{ {\mathcal{F}} }

\newcommand{\cI}{ {\mathcal{I}} }
\newcommand{\cJ}{ {\mathcal{J}} }
\newcommand{\cL}{ {\mathcal{L}} }

\newcommand{\cP}{ {\mathcal{P}} }

\newcommand{\cR}{ {\mathcal{R}} }

\newcommand{\anorme}{ {\mathopen{}\left\lVert \cdot \right\rVert\mathclose{}} }
\newcommand{\norme}[1]{ {\mathopen{}\left\lVert #1 \right\rVert\mathclose{}} }

\newcommand{\aabs}{ {\mathopen{}\left|\cdot\right|} }

\newcommand{\indi}[1]{ {\mathbbm{1}_{#1}} }

\newcommand{\iense}[2]{ {\llbracket #1 , #2 \rrbracket } }
\newcommand{\ense}[1]{ \set*{#1} } 													
\newcommand{\inteoo}[2]{ { \mathopen{}\left( #1  ,  #2  \right)\mathclose{} } }

\newcommand{\inteff}[2]{ { \mathopen{}\left[ #1  ,  #2  \right]\mathclose{} } }

\newcommand{\ade}[1]{ {\overline{#1}} }
\newcommand{\fr}[1]{ {\partial #1} }

\newcommand{\regu}[2]{ {\cC^{#1}\lrpa{#2}}  }
\newcommand{\reguc}[2]{ {\cC_0^{#1}\lrpa{#2}}  }
\newcommand{\regub}[2]{ {\cC_b^{#1}\lrpa{#2}}  }

\newcommand{\td}{ {\text{d}} }
\newcommand{\dif}[3]{ {\text{d}^{#3}_{#2}\, {#1}} }
\newcommand{\Dif}[3]{ {\text{D}^{#3}_{#2}\, {#1}} }

\newcommand{\devp}[2]{ { \lrbr{\partial_{#2} #1} } }

\newcommand{\Integd}[4]{ {\int_{#2}^{#3} #1 \,\text{d}#4} }

\newcommand{\integd}[3]{ {\int_{#2} #1 \,\text{d}#3} }
\newcommand{\integdd}[4]{ {\int_{#2} #1 \,\text{d}#3\text{d}#4} }
\newcommand{\integ}[2]{ {\int_{#2} #1} }

\newcommand{\dive}{ {\text{div}} }

\newcommand{\supp}{ {\text{supp}} }

\renewcommand{\devp}[2]{ { \partial_{#2} #1 } }

\numberwithin{equation}{section}

\tcbset{colback=cyan!5!white,colframe=cyan!75!black}

\newcommand{\zH}{ {\prescript{}{0}{H}} }
\newcommand{\tH}{ {\prescript{}{T}{H}} }

\newcommand{\hdirz}{ {\zH^{\frac{3}{4},\frac{3}{2}}(\gS)} }
\newcommand{\hdirt}{ {\tH^{\frac{3}{4},\frac{3}{2}}(\gS)} }

\newcommand{\gvpSpace}{ G }
\newcommand{\ASpace}{ \cA }
\newcommand{\SolSpace}{ H^{1,2}(Q) }
\newcommand{\qSpace}{ L^\infty(Q)}

\title{\textbf{Recovery of coefficients for a convection-diffusion equation from partial data}}

\author{Liam Buisson\footnote{Univ Rouen Normandie, Normandie Univ, LMRS UMR 6085, F-76000, Rouen, France; \\ \indent \indent E-mail: liam.buisson1@univ-rouen.fr}}
\date{}

\begin{document}


\maketitle

\begin{abstract}
	{\scriptsize This article is devoted to the inverse problem of determining the zeroth- and first-order coefficients, depending on both the time and space variables, in a parabolic equation from partial boundary measurements of the flux generated by Dirichlet excitations. 
	More precisely, we establish the unique determination of a time-dependent convection term and potential from the partial Dirichlet-to-Neumann map associated with the corresponding parabolic equation, where the Neumann measurements are restricted to the portion of the boundary illuminated from a point located outside the closure of the domain. 
	Our main objective is to extend to parabolic equations several observations that have thus far been confined to the elliptic setting and to exploit these properties to recover a general class of coefficients depending on both time and space variables. 
    To achieve this, we develop a suitable class of special solutions to the parabolic equation, introduce a new family of Carleman estimates with nonlinear weights, and combine these tools with partial data results for the X-ray transform. 
	The principal difficulty of the present work stems from the fact that, unlike the existing literature on partial data inverse problems for parabolic equations, both the phase of the special solutions and the weight of the corresponding Carleman estimates are nonlinear. 
	As a byproduct of our analysis, we establish a new class of Carleman estimates for parabolic equations with nonlinear weights, which may be viewed as the parabolic analogue of the notion of limiting Carleman weights in the elliptic setting.
	
	\noindent
	{\bf  Keywords:} Inverse problem, parabolic equations, convection-diffusion equation, uniqueness, Carleman estimates, partial data.

	\noindent
	{\bf Mathematics subject classification 2020 :} 35R30, 35K20.
	}

    \noindent
\end{abstract}

\tableofcontents


\section{Introduction:}

Let $\gO \subset \RR^n$ ($n\geq 2$) be a bounded open set such that $\RR^n \setminus \gO$ is connected and whose boundary $\fr{\gO} =: \gG$ is of class $\cC^2$.
We fix $T>0$ and set $Q=\inteoo0T\times \gO$ and $\gS=\inteoo0T\times \gG$. 
Then, we consider the following initial-boundary value problem (IBVP):
\begin{equation}
    \label{Eq:P}
    \begin{dcases}
        \partial_t u(t,x)-\gD u(t,x)+A(t,x)\cdot\nabla u(t,x)+q(t,x)u(t,x)=0,
        & (t,x)\in Q,\\
        u(t,s)=g(t,s),
        & (t,s)\in \gS,\\
        u(0,x)=0,
        & x\in\gO.
    \end{dcases}
\end{equation}

For every $y\in\RR^n\setminus\ade{\gO}$, we define the shadowed and illuminated portions of the boundary $\gG$, as viewed from the point $y$ defined by
$$ \gG_\pm(y):= \ense{ x\in\gG \given \pm (x-y)\cdot\nu(x)\geq 0 }, $$
where $\nu(x)$ denotes the unit outward normal vector to $\gG$ at $x$, and $\cdot$ (resp. $\aabs$) stands for the Euclidean inner product (resp.\ norm) in $\RR^n$.

Let $x_0\in\RR^n\setminus\overline{\gO}$ and let $U$ be a neighborhood of $\gG_-(x_0)$ in $\gG$. 
We formally define the partial parabolic Dirichlet-to-Neumann map (DN map) associated with \eqref{Eq:P} by
\begin{equation}
    \label{Eq:pdn}
    \gL_{(A,q)} : g\longmapsto \devp{u}{\nu}_{|_{\inteoo0T \times U}},
\end{equation}
where $u$ is the unique solution to \eqref{Eq:P}. 
We refer to Corollary~\ref{Cor:WP1} for a rigorous definition of the partial DN map \eqref{Eq:pdn}.

In the present article, we study the following inverse problem:
\begin{itemize}
\item[{\bf(IP)}]
\emph{Determine uniquely and simultaneously the coefficients $A$ and $q$ from the knowledge of the partial DN map $\gL_{(A,q)}$.}
\end{itemize}

\subsection{Motivations:}

Second-order parabolic equations of the form \eqref{Eq:P} belong to the class of convection--diffusion--reaction equations (also referred to as advection--diffusion--reaction equations). 
They arise in the modeling of a wide variety of physical phenomena, including bond pricing (through the Black--Scholes equation), the time evolution of probability densities associated with diffusion processes (e.g., the Fokker--Planck equation \cite{riskenFokkerPlanckEquationMethods1996}), and the evolution of the density of a physical quantity in systems such as nutrient transport in ecosystems, population dispersal in heterogeneous environments, or the dispersion of atmospheric pollutants (see, for instance, \cite{stockerIntroductionClimateModelling2011,okuboDiffusionEcologicalProblems2001,seinfeldAtmosphericChemistryPhysics1986}).
In these models, the solution $u$ represents the transported quantity (e.g., concentration in mass transfer or temperature in heat transfer), the vector field $A$ describes the underlying velocity field, and the term $-qu$ models the local generation or depletion rate of $u$.
Moreover, in many practical applications, both $A$ and $q$ depend on time and space.

From this perspective, inverse problem {\bf(IP)} consists in recovering the velocity field $A$ and the reaction coefficient $q$ from measurements of the boundary flux on a neighborhood of the portion $\gG_-(x_0)$ of the boundary illuminated  by $x_0$, generated by suitable boundary excitations.

\subsection{State of the art:}

Let us first remark that the inverse problem \textbf{(IP)} belongs to the class of inverse coefficient problems, which are distinguished both by their broad range of applications and by their intrinsic mathematical interest stemming from their nonlinearity and ill-posedness.
One of the pioneering contributions in this direction is due to \cite{sylvesterGlobalUniquenessTheorem1987}, who gave the first positive answer to the celebrated Calderón problem. 
Since then, inverse coefficient problems have been extensively investigated, leading to numerous uniqueness and stability results for the recovery of first- and zeroth-order coefficients in elliptic equations (see, e.g., \cite{krupchykUniquenessInverseBoundary2014,krupchykInverseProblemsMagnetic2018,NSU,sunInverseBoundaryValue1993}).

A central topic in inverse coefficient problems is the notion of partial data, where the objective is to relax the requirement of full boundary measurements by restricting either the set of boundary excitations, the set of observations, or both, to suitable portions of the boundary.
In problem \textbf{(IP)}, the partial data setting arises from the fact that the Neumann measurements are available only on the subset $U$ of the boundary $\gG$.
This restriction is motivated by the desire to reduce the large amount of information encoded in the full parabolic Dirichlet-to-Neumann map (i.e. the map $\gL_{(A,q)}$ when $U=\gG$). 
For elliptic equations, partial data inverse problems have been the subject of intensive research.
One of the first contributions in this direction is due to \cite{bukhgeimRecoveringPotentialPartial2002}, whose analysis was substantially extended in \cite{kenigCalderonProblemPartial2007} to the recovery of zeroth-order coefficients under partial boundary measurements similar to those considered in \textbf{(IP)}. 
The methods developed in \cite{bukhgeimRecoveringPotentialPartial2002,kenigCalderonProblemPartial2007} were subsequently generalized to the simultaneous recovery of first- and zeroth-order coefficients in several works, including \cite{Ch,DKSU,LRT,Sel}.

For parabolic equations, most of the available partial data results concern time-independent coefficients, for which the problem can be reduced either to inverse spectral problems or to inverse problems for a family of elliptic equations (see, e.g., \cite{CaKa,KKLM,KiIP,KOSY,KLLY}).
Such reduction techniques are no longer applicable in the presence of time-dependent coefficients.
Nevertheless, several alternative approaches have been developed to address inverse coefficient problems of the type \textbf{(IP)} for parabolic equations with time-dependent coefficients.
One of the pioneering works in this direction is due to \cite{isakovCompletenessProductsSolutions1991}, who established a density property for products of solutions to parabolic equations and used it to prove uniqueness for the recovery of a time-dependent zeroth-order coefficient.
This result was significantly strengthened by \cite{choulliLogarithmicStabilityDetermining2018}, who established both stability and partial data results by constructing suitable solutions using Carleman estimates with linear weights, following the strategy introduced for elliptic equations in \cite{bukhgeimRecoveringPotentialPartial2002,kenigCalderonProblemPartial2007}.
The approach of \cite{choulliLogarithmicStabilityDetermining2018} was further extended by \cite{caroDeterminationConvectionTerms2019} to the simultaneous recovery of first- and zeroth-order coefficients with low regularity assumptions.
Since then, the methodology developed in \cite{choulliLogarithmicStabilityDetermining2018,caroDeterminationConvectionTerms2019} has been applied in several works devoted to the recovery of time-dependent coefficients (see, e.g., \cite{bellassouedStabilityEstimateInverse2020,SV,Ra}).
We also mention the recent work \cite{feizmohammadiInverseBoundaryValue2023}, where the author investigated the recovery of leading-order time-dependent coefficients by constructing a new class of special solutions for parabolic equations with variable principal coefficients.
This analysis also provides a new construction of solutions for parabolic equations that was subsequently exploited by \cite{kianDeterminingBreakingGauge2024} to recover a time-dependent zeroth-order coefficient in dimension $n \geq 3$ using the data associated with problem \textbf{(IP)}.
Beyond their intrinsic interest, inverse problems involving time-dependent coefficients also play a fundamental role in the study of inverse problems for nonlinear equations (see, e.g., \cite{FKU,Isa93,kianDeterminingBreakingGauge2024,KiU}).

The main objective of the present article is to extend to parabolic equations the partial data framework developed for elliptic equations in \cite{DKSU,kenigCalderonProblemPartial2007}. 
More precisely, we consider measurements of the Neumann data restricted to a neighborhood of the  part of the boundary illuminated by $x_0 \in \mathbb{R}^n \setminus \ade{\gO}$, as described in problem \textbf{(IP)}. 
Our approach relies on extending the construction of special solutions introduced in \cite{feizmohammadiInverseBoundaryValue2023} to convection--diffusion equations of the form \eqref{Eq:P}, together with the derivation of a new family of Carleman estimates for parabolic equations with variable first-order coefficients stated in Appendix~\ref{Carleman}.
In particular, we establish a parabolic analogue of the notion of limiting Carleman weights introduced in \cite{dossantosferreiraLimitingCarlemanWeights2009} for nonlinear weights.
We believe that this new class of Carleman estimates with nonlinear weights is of independent interest and may lead to further applications, including the simultaneous restriction of both the observation set and the Dirichlet excitation set, in the spirit of \cite{kenigCalderonProblemPartial2007}.

\subsection{Outline:}

This article is organized as follows.
In Section~\ref{MaRe}, we state and comment the main results of this article.
We recall in Section~\ref{WPIBVP}, the well-posedness of \eqref{Eq:P} in $H^1\lrpa{ 0,T; L^2(\gO) } \cap L^2\lrpa{ 0,T; H^2(\gO) }$, for bounded coefficients $A$ and $q$ and suitable regularity of $g$.
In Section~\ref{Gauge}, we show the invariance of the partial DN map \eqref{Eq:pdn} for a certain gauge transformation.
In Section~\ref{GOS}, we construct particular solutions of problem \eqref{Eq:P} and its adjoint.
Such solutions will be used in Section~\ref{UnPr} to prove Theorem~\ref{Theo1}.
After showing that the knowledge of the partial DN map \eqref{Eq:pdn} determines the gauge class of $(A,q)$ (see Section~\ref{Gauge} for definitions), we give the uniqueness of this inverse problem under further assumptions on $A$ and $q$ in Section~\ref{FuCo}.
The appendix of this article is devoted to the proof of useful lemmas for our main results, along with a Carleman estimate for a type of weight function that we will specify both in Section~\ref{GOS} and at the beginning of Appendix~\ref{Carleman}.

\section{Main results:}			\label{MaRe}

Let us recall that there is a natural obstruction to inverse problem \textbf{(IP)}, already observed by \cite{caroDeterminationConvectionTerms2019}. 
More precisely, problem \textbf{(IP)} enjoys a gauge invariance that can be formally described as follows. For any function $\gvp$, belonging to a suitable function space, satisfying
$\gvp_{|_\gS}=0$ and $\devp{\gvp}{\nu}_{|_\gS}=0$,
the Dirichlet-to-Neumann map is invariant under the gauge transformation
$$ S_\gvp:(A,q)\longmapsto \left( A+2\nabla\gvp,\, q-\left(\devp{\gvp}{t}-\gD\gvp+|\nabla\gvp|^2+A\cdot\nabla\gvp\right) \right). $$
We refer to Section~\ref{Gauge} for a rigorous statement and proof of this invariance. 
Throughout the paper, we denote by the \emph{gauge class} of $(A,q)$ the set of all pairs $S_\gvp(A,q)$, where $\gvp$ belongs to the function space $G$ introduced in Section~\ref{Gauge}.

In view of this obstruction, problem \textbf{(IP)} naturally splits into the following two inverse problems:
\begin{itemize}
    \item[{\bf(IP1)}] \emph{Determine the gauge class of $(A,q)$ from the knowledge of $\gL_{(A,q)}$.}
    \item[{\bf(IP2)}] \emph{Under which additional assumptions on $(A,q)$ can one uniquely recover $(A,q)$ from the knowledge of $\gL_{(A,q)}$?}
\end{itemize}

We study these two problems separately, beginning with \textbf{(IP1)}.
For every $y\in \RR^n\setminus\ade{\gO}$ and $\eta \geq 0$, we define
\[
\gG_\pm(y,\eta) := \ense{ x\in\gG \given \pm(x-y)\cdot\nu(x)>\pm\eta }, \qquad \gS_\pm(y,\eta) := \inteoo0T \times\gG_\pm(y,\eta).
\]
From now on, without loss of generality, we fix $U=\gG_-(x_0,2\gve)$ in the definition \eqref{Eq:pdn} of $\gL_{(A,q)}$. 
We assume that $\gve>0$ is sufficiently small so that
\begin{equation}
\label{Hyp1}
B(x_0,\gve)\cap \lrpa{\gO+B(0,\gve)}=\emptyset,
\end{equation}
where $B(a,\eta)$ denotes the open Euclidean ball of center $a$ and radius $\eta$. 
We also introduce the space
\[
\ASpace
:=
\left(
\regub{0}{0,T;\regu{2}{\ade{\gO}}}
\cap
\regub{1}{0,T;\regu{0}{\ade{\gO}}}
\right)^n,
\]
whose precise definition is recalled in Section~\ref{Gauge}. 
Our answer to problem \textbf{(IP1)} is the following.

\begin{Theoreme}
\label{Theo1}
Let $A_1,A_2\in\ASpace$ and $q_1,q_2\in L^\infty(Q)$. If
\[
\gL_{(A_1,q_1)}
=
\gL_{(A_2,q_2)}
\quad\text{and}\quad
{A_1}_{|_\gS}
=
{A_2}_{|_\gS},
\]
then there exists $\gvp\in\gvpSpace$ such that
\[
(A_2,q_2)=S_\gvp(A_1,q_1).
\]
\end{Theoreme}

To the best of our knowledge, Theorem~\ref{Theo1} provides the first positive answer to problem \textbf{(IP1)} with Neumann measurements restricted to a neighborhood of the illuminated part of the boundary associated with an arbitrary point
$x_0\in \RR^n\setminus\ade{\gO}$, namely the same measurement configuration considered for elliptic equations in \cite{DKSU,kenigCalderonProblemPartial2007}. 
We recall that, when $\ade{\gO}$ is strictly convex, the set $U$ appearing in \eqref{Eq:pdn} may be chosen as an arbitrary open subset of $\gG$, making this class of restrictions optimal. 
Even in the case $A\equiv0$, Theorem~\ref{Theo1} improves the existing literature by covering the two-dimensional case, which was not addressed in \cite[Corollary~2.4]{kianDeterminingBreakingGauge2024}.

Our proof relies on the combination of a suitable class of solutions to \eqref{Eq:P}, usually referred to as geometric optics (GO) solutions, constructed in Section~\ref{GOS}, together with a new family of Carleman estimates established in Appendix~\ref{Carleman}. 
The main difficulty in the analysis of \textbf{(IP1)} is that, unlike the existing approaches of \cite{choulliLogarithmicStabilityDetermining2018,caroDeterminationConvectionTerms2019}, our construction requires Carleman estimates with nonlinear weights appearing in the phase of the GO solutions. 
To overcome this difficulty, we extend to convection--diffusion equations of the form \eqref{Eq:P} the energy-based methodology introduced in \cite{feizmohammadiInverseBoundaryValue2023}. 
This is combined with a new class of Carleman estimates adapted to nonlinear weights. 
As a byproduct, we establish in Appendix~\ref{Carleman} a parabolic analogue of the notion of limiting Carleman weights introduced in \cite{dossantosferreiraLimitingCarlemanWeights2009}. 
We believe that this result is of independent interest and may have further applications, including the simultaneous restriction of both the observation set and the Dirichlet excitation set, in the spirit of \cite{kenigCalderonProblemPartial2007}.
Our analysis also relies on a suitable approximation of the unknown coefficients, enabling us to relax the assumption that the first-order derivatives of coefficients $A_j$, $j=1,2$, are known on the lateral boundary $\gS$. 
To the best of our knowledge, this approach is novel and is presented in Section~\ref{Flem}.
\new
 
We now turn to problem \textbf{(IP2)}.
Since $\gO$ is sufficiently smooth, for every $t\in \inteoo{0}{T}$ the Hodge decomposition of $A(t,\cdot)$ is unique and takes the form
$$ A(t,\cdot)=\nabla\phi+R \quad\text{in }\gO, $$
where $\phi\in H^2(\gO)$ satisfies $\phi_{|_\gG}=0$ and $R\in H^1(\gO)^n$ satisfies $\dive(R) =0$ in $\gO$.
Under this decomposition, Theorem~\ref{Theo1} shows that the partial Dirichlet-to-Neumann map \eqref{Eq:pdn} determines the solenoidal part of $A$, provided that $A$ is sufficiently smooth. 
However, because of the gauge invariance described above, the gradient component cannot be recovered, and therefore neither can the pair $(A,q)$.

This observation provides a first answer to \textbf{(IP2)}. 
Indeed, if
\[
(A',q')=S_\gvp(A,q)
\]
for some $\gvp\in\gvpSpace$, then for every $t\in \inteoo0T$ the function
$\phi=\gvp(t,\cdot)$ satisfies
$$
\begin{dcases}
    \gD\phi(x) = \frac12 \dive \bpa{(A'-A)(t,\cdot)}(x),    &x\in\gO,\\
    \phi(x)=0,  &x\in\gG.
\end{dcases}
$$
Therefore, if $A'-A$ is divergence free in $Q$, then necessarily $\gvp=0$, which implies $(A,q)=(A',q')$. 
Although this criterion is based on the relation between $A$ and $A'$, we show in Section~\ref{FuCo}, through Corollaries~\ref{Coro1} and \ref{Coro2}, that suitable relations between $q$ and $q'$ may also force $\gvp$ to vanish by means of a unique continuation argument (see Lemma~\ref{Lem1}).

Our second main result provides another answer to problem \textbf{(IP2)}.

\begin{Theoreme}
\label{Theo2}
For every $A\in\ASpace$ and every $q_1,q_2\in\qSpace$,
\[
\gL_{(A,q_1)}
=
\gL_{(A,q_2)}
\quad\Longrightarrow\quad
q_1=q_2.
\]
\end{Theoreme}

As shown in Section~\ref{UnPr}, Theorem~\ref{Theo2} also constitutes an intermediate step in the proof of Theorem~\ref{Theo1}. 
Its proof relies on the geometric optics solutions constructed in Section~\ref{GOS} together with uniqueness results for the X-ray transform of compactly supported distributions (see \cite{ilmavirtaUniqueContinuationNormal2020}).

When the zeroth-order coefficient is known, we obtain the following consequence.

\begin{Corollaire}
\label{Coro1}
Let $q\in L^\infty(Q)$ and let
$A_1,A_2\in\ASpace$ satisfy
${A_1}_{|_\gS}={A_2}_{|_\gS}$.
Then
\[
\gL_{(A_1,q)}
=
\gL_{(A_2,q)}
\quad\Longrightarrow\quad
A_1=A_2.
\]
\end{Corollaire}

Our final result concerns the Fokker--Planck equation
\begin{equation}
\label{PFK1}
\begin{dcases}
\devp{u}{t}-\gD u+\dive(Au)=0,
&\text{in }Q,\\
u_{|_\gS}=g,
&\text{on }\gS,\\
u_{|_{t=0}}=0,
&\text{in }\gO,
\end{dcases}
\end{equation}
together with the boundary operator
\[ \cF_A:g\longmapsto \devp{u}{\nu}_{|_{\gS_-(x_0,2\gve)}}, \]
where $u$ denotes the unique solution of \eqref{PFK1}.

\begin{Corollaire}
\label{Coro2}
Let $A_1,A_2\in\ASpace$ satisfy
${A_1}_{|_\gS}={A_2}_{|_\gS}$.
Then
\[ \cF_{A_1} = \cF_{A_2} \quad\Longrightarrow\quad A_1=A_2. \]
\end{Corollaire}

To the best of our knowledge, Theorem~\ref{Theo2} together with Corollaries~\ref{Coro1} and \ref{Coro2} provides the first positive answer to problem \textbf{(IP2)} for parabolic equations with partial boundary measurements of the type considered in \cite{DKSU,kenigCalderonProblemPartial2007}. 
These results naturally complement Theorem~\ref{Theo1} by showing that, under suitable additional assumptions, the gauge invariance can be broken, leading to the complete recovery of either the convection field $A$ or the zeroth-order coefficient $q$. 
Such results are particularly relevant in applications where only one of these two parameters is of interest.

\section{Well-posedness of the IBVP:}      \label{WPIBVP}

In this section, we recall some properties of problem \eqref{Eq:P} and we define the partial DN map \eqref{Eq:pdn}.
We start with some preliminary tools.
For $r \geq 0$, $s \in \inteff{0}{2}$ and $X$ being either $\gO$, $\gG$ or $\gG_-(x_0,2\gve)$, we define the real Hilbert space $H^{r,s}\lrpa{ \inteoo{0}{T} \times X }$ as
$$ H^{r,s}\lrpa{ \inteoo{0}{T} \times X } := H^r\bpa{ 0,T; L^2(X) } \cap L^2\bpa{ 0,T; H^s(X) }, $$
with the following norm induced by its inner product 
$$ \forall u \in H^{r,s}\lrpa{ \inteoo{0}{T} \times X } ,\quad \norme{u}_{H^{r,s}\lrpa{ \inteoo{0}{T} \times X }}^2 := \norme{u}_{H^r (0,T; L^2(X))}^2 + \norme{u}_{L^2(0,T; H^s(X))}^2. $$
We now write $H^{r,s}(Q)$ (resp. $H^{r,s}(\gS)$) (resp. $H^{r,s}(\gS_-(x_0,2\gve))$) instead of $H^{r,s}\lrpa{ \inteoo{0}{T} \times \gO }$ (resp. $H^{r,s}\lrpa{ \inteoo{0}{T} \times \gG}$) (resp. $H^{r,s}(\inteoo{0}{T}\times\gG_-(x_0,2\gve))$). 
In addition, we set
$$ \prescript{}{t_0}{H}^{r,s}(\gS) := \ense{ u \in H^{r,s}(\gS) \given u_{|_{t=t_0}} = 0 }, $$
when $t_0$ is either $0$ or $T$, and $r > \frac{1}{2}$.
Let us consider the following proposition,

\begin{Proposition}
    \label{Pro:WP1}
    For fixed $A \in L^\infty(Q)^n$, $q \in L^\infty(Q)$, $g \in \hdirz$ and $F \in L^2(Q)$, the following IBVP
    \begin{equation}
        \label{Pl}
        \begin{dcases}
            \devp{u}{t} - \gD u + A\cdot\nabla u + qu = F, &\quad \text{ in } Q,\\
            u_{|_\gS} = g, &\quad \text{ on } \gS, \\
            u_{|_{t=0}} = 0, &\quad \text{ on } \gO,
        \end{dcases}
    \end{equation}
    admits a unique solution $u_{g,F} \in H^{1,2}(Q)$, and the linear map 
    $$ S :  \hdirz \times L^2(Q) \ni (g,F) \longmapsto u_{g,F} \in H^{1,2}(Q)$$
    is continuous.
\end{Proposition}

Applying \cite[Chapter 4, Theorem 2.1]{lionsNonHomogeneousBoundaryValue1972a}, we deduce the following corollary,
\begin{Corollaire}
    \label{Cor:WP1}
    For every $A \in L^\infty(Q)^n$ and $q \in L^\infty(Q)$, the map $\gL_{(A,q)}$ in \eqref{Eq:pdn} is well-defined and continuous from $\hdirz$ to $H^{\frac14,\frac12}(\gS_-(x_0,2\gve)) $.
\end{Corollaire}
Another consequence of Proposition~\ref{Pro:WP1} is the well-posedness for the formal adjoint problem.
\begin{Corollaire}
    \label{Cor:WP2}
    For fixed $A \in L^\infty(Q)^n$, $q \in L^\infty(Q)$, $g \in \hdirt$ and $F \in L^2(Q)$ the following final boundary value problem
    \begin{equation}
        \begin{dcases}
            - \devp{v}{t} - \gD v + A\cdot\nabla v + qv = F, &\quad \text{ in } Q,\\
            v_{|_\gS} = g, &\quad \text{ on } \gS, \\
            v_{|_{t=T}} = 0, &\quad \text{ on } \gO,
        \end{dcases}
    \end{equation}
    admits a unique solution $v_{g,F} \in H^{1,2}(Q)$, and the linear map 
    $$ \hdirt \times L^2(Q) \ni (g,F) \longmapsto v_{g,F} \in H^{1,2}(Q)$$
    is continuous.
\end{Corollaire}
\new

We will now consider the proof of Proposition~\ref{Pro:WP1}. 
This result is classical but we give its proof for sake of completeness.
\begin{Preuve}[of Proposition~\ref{Pro:WP1}]
    Let us fix $A \in L^\infty(Q)^n$ and $q \in L^\infty(Q)$ and consider $F \in L^2(Q)$. 
    We first demonstrate that the problem \eqref{Pl} is well-posed in the case of a homogeneous boundary condition.
    Then, we conclude by using lifting argument.
    
    \textbf{Step 1 - Homogeneous boundary condition:}
    We denote by $(\cdot|\cdot)$ the duality bracket over $ L^2\lrpa{0,T;H^1_0(\gO)}$.
    Applying \cite[Chapter 3, Theorem 4.1]{lionsNonHomogeneousBoundaryValue1972}, for every $f$ in $L^2\lrpa{0,T;H^{-1}(\gO)}$ (which we identify to $L^2\lrpa{0,T;H^1_0(\gO)}'$), the variational problem 
    \begin{equation}
        \label{IB1}
        \begin{dcases}
            (\devp{u}{t}|v) + \integdd{\nabla u \cdot \nabla v + A\cdot \nabla u v + quv}{Q}{x}{t} = (f | v), \quad \forall v \in L^2\lrpa{0,T;H^1_0(\gO)}, \\
            u_{|_{t=0}} = 0,
        \end{dcases}
    \end{equation}
    admits a unique solution $u_f$ in the Hilbert space $W(0,T) := L^2\lrpa{0,T;H^1_0(\gO)} \cap H^1\lrpa{0,T;H^{-1}(\gO)}$ and the map $f \mapsto u_f$ is continuous from $L^2\lrpa{0,T;H^{-1}(\gO)}$ onto $W(0,T)$ when the latest space is equipped with the following norm,
    $$ \forall u \in W(0,T), \quad \norme{u}_{W(0,T)}^2 = \norme{u}^2_{L^2\npa{0,T;H^1_0(\gO)}} + \norme{\devp{u}{t}}^2_{L^2\npa{0,T;H^{-1}(\gO)}} .$$
    If we consider $f := F - A \cdot \nabla u_F - qu_F \in L^2(Q)$, then \cite[Chapter 4, Theorem 1.1]{lionsNonHomogeneousBoundaryValue1972a} gives us the existence and uniqueness of $u \in H^{1,2}(Q) \cap L^2(0,T;H^1_0(\gO))$ being solution of the following IBVP
    \begin{equation*}
        \begin{dcases}
            \devp{u}{t}- \gD u = F - A \cdot \nabla u_F - qu_F , &\quad \text{ in } Q,\\
            u_{|_\gS} = 0, &\quad \text{ on } \gS, \\
            u_{|_{t=0}} = 0, &\quad \text{ on } \gO.
        \end{dcases}
    \end{equation*}
    Since $F$ is also in $L^2\lrpa{0,T;H^{-1}(\gO)}$, $u$ is a solution of the variational problem
    \begin{equation*}
        \begin{dcases}
            (\devp{u}{t}|v) + \integdd{\nabla u \cdot \nabla v}{Q}{x}{t} = (F | v) - \integdd{\lrpa{A\cdot \nabla u_F + qu_F}v}{Q}{x}{t}, \quad \forall v \in L^2\lrpa{0,T;H^1_0(\gO)}, \\
            \tilde u_{|_{t=0}} = 0,
        \end{dcases}
    \end{equation*}
    whose unique solution is $u_F$. 
    We then deduce that $u_F \in H^{1,2}(Q)\cap L^2\lrpa{0,T;H^1_0(\gO)}$ and the problem
    \begin{equation}
        \label{Pll}
        \begin{dcases}
            \devp{u}{t} - \gD u + A \cdot \nabla u + qu = F , &\quad \text{ in } Q,\\
            u_{|_\gS} = 0, &\quad \text{ on } \gS, \\
            u_{|_{t=0}} = 0, &\quad \text{ on } \gO,
        \end{dcases}
    \end{equation}
    admits a unique solution $u_F \in H^{1,2}(Q)\cap L^2\lrpa{0,T;H^1_0(\gO)}$ for every $F \in L^2(Q)$.
    Moreover, by construction and, following \cite[Chapter 4, Theorem 1.1 and Remark 1.6]{lionsNonHomogeneousBoundaryValue1972a}, we know that there exists $C \geq 0$ such that for every $F \in L^2(Q)$, we have
    $$ \norme{u_F}_{H^{1,2}(Q)} \leq C \norme{F - A\cdot \nabla u_F - qu_F}_{L^2(Q)} .$$
    We then obtain, thanks to the well-posedness of \eqref{IB1} from $L^2(Q)$ to $L^2(0,T;H^1_0(\gO))$,
    \begin{align*}
        \norme{u_F}_{H^{1,2}(Q)} &\leq C \norme{F - A\cdot \nabla u_F - qu_F}_{L^2(Q)} \\
        &\leq C\lrbr{\norme{F}_{L^2(Q)} + \lrpa{\norme{A}_{L^\infty(Q)^n} + \norme{q}_{L^\infty(Q)}}\norme{u_F}_{L^2\lrpa{0,T;H^1_0(\gO)}}} \\
        &\leq C\lrbr{1+\lrpa{\norme{A}_{L^\infty(Q)^n} + \norme{q}_{L^\infty(Q)}}C'}\norme{F}_{L^2(Q)},
    \end{align*}
    for some constant $C' \geq 0$ depending only on $\gO$ and $T$.
    \new

    \textbf{Step 2 - Non-homogeneous boundary condition:}
    The uniqueness in the homogeneous boundary condition case and the linearity of \eqref{Pl} implies the uniqueness in the case of a non-homogeneous boundary condition.
    Applying \cite[Chapter 4, Theorem 2.3]{lionsNonHomogeneousBoundaryValue1972a}, the linear map 
    \begin{equation}
        \label{IB2}
        \zH^{1,2}(Q) \ni u \longmapsto u_{|_\gS} \in \hdirz
    \end{equation}
    is continuous and surjective. 
    Since $\zH^{1,2}(Q)$ is a Hilbert space, there exists a continuous linear map $\cR : \hdirz \rightarrow \zH^{1,2}(Q) $ such that,
    $$ \forall g \in \hdirz,\quad \cR(g)_{|_\gS} = g. $$
    Using this, we consider $g\in \hdirz$ and $F \in L^2(Q)$ and set $v = \cR(g)$ along with 
    $$f := F - (\devp{v}{t} -\gD v + A\cdot\nabla v + qv) \in L^2(Q).$$
    Then, one can check that $u_{F,g} := u_f + \cR(g) \in H^{1,2}(Q)$ solves \eqref{Pl}.
    This proves the existence of a solution to \eqref{Pl}.
    The continuity of $S$ follows from the continuity of $\cR$, the continuity of 
    $$ \hdirz \times L^2(Q) \ni (g,F) \longmapsto F - (\devp{}{t} -\gD + A\cdot\nabla + q)\cR(g) \in L^2(Q), $$
    and the well-posedness of \eqref{Pll}.
\end{Preuve}

\section{The gauge invariance:}            \label{Gauge}

In this section, we detail the fact that we cannot determine $(A,q)$ from the sole knowledge of $\gL_{(A,q)}$. 
We will also define more explicitly the notion of gauge transformation, and we will set the framework of this article.
We first recall some preliminary functional spaces.
For any non-empty open subset $U$ of $\RR^n$, we define for $k = 1,2$,
$$\regub{k}{U} := \ense{ \gvp \in \regu{k}{U} \given \forall \ga \in (\NN\cup \ense{0})^n, \quad \abs{\ga} \leq k, \quad\sup_U \abs{\devp{^\ga \gvp}{}} < \infty}, $$ 
that we consider with its norm $\anorme_{\regub{k}{U}}$ defined by $\norme{\gvp}_{\regub{k}{U}} := \sum_{\abs{\ga} \leq k} \sup_U \abs{\devp{^\ga \gvp}{}}$ for every $\gvp$ in $\regu{k}{U}$. 
We denote by $\regu{1}{\ade{\gO}}$ the space of all functions which are restrictions to $\gO$ of a function in $\regub{1}{\RR^n}$, that we consider with the following norm 
$$\norme{\gvp}_{\regu{1}{\ade{\gO}}} = \inf\ense{ \norme{\tilde{\gvp}}_{\regub{1}{\RR^n}} \given \forall \tilde{\gvp} \in \regub{1}{\RR^n}, \quad \tilde{\gvp}_{|_\gO} = \gvp}. $$
Notice that $\gvp_{|_\gS}$ and $(\nabla \gvp)_{|_\gS}$ are well-defined over $\gS$ for every $\gvp \in \regu{1}{\ade{\gO}}$, using a continuity argument.
We consider $\regu{1}{\ade{\gO}} \cap \regub{2}{\gO}$ with $\anorme_{\regu{1}{\ade{\gO}} \cap \regub{2}{\gO}} := \anorme_{\regu{1}{\ade{\gO}}} + \anorme_{\regub{2}{\gO}}$.
We point out the fact that, since $\gO$ have a $\cC^1$ boundary, we know that $\regub{2}{\gO}$ is continuously embedded into $\regu{1}{\ade{\gO}}$. 
We will also consider the space $\regu{0}{\RR^n}$ and $\regu{2}{\RR^n}$ with their usual locally convex space structure.
For a locally convex space $\lrpa{E,\ense{p_i}_{i\in I}}$, we define 
$$ \regub{0}{0,T;E} := \ense{ f \in \regu{0}{0,T;E} \given \forall i \in I,\quad \exists M_i \geq 0, \quad \sup_{t\in \inteoo{0}{T}} p_i(f(t)) \leq M_i}. $$
We also define the following subspace of $\regub{0}{0,T;E}$,
$$ \regub{1}{0,T;E} := \ense{ f \in \regu{1}{0,T;E} \given \devp{f}{t} \in \regub{0}{0,T;E} }. $$
\new

Let us fix $A \in L^\infty(Q)^n$ and $q \in L^\infty(Q)$, the set $C := \regub{1}{0,T;\regu{0}{\ade{\gO}}} \cap \regub{0}{0,T;\regu{1}{\ade{\gO}} \cap \regub{2}{\gO}}$, along with the gauge function space:
$$\gvpSpace := \ense{ \gvp \in C \given \gvp_{|_\gS} = 0 \text{ and } \devp{\gvp}{\nu}_{|_\gS} = 0}. $$
We also define the operator $L_{(A,q)} := \devp{t}{} - \gD + A\cdot \nabla + q$.
For $\gvp \in \gvpSpace$ and $v \in H^{1,2}(Q)$, fixing $u = e^{-\gvp}v$ gives $u_{|_\gS} = v_{|_\gS}$, $\devp{u}{\nu}_{|_\gS} = \devp{v}{\nu}_{|_\gS}$ and
$$ L_{(A,q)}u = 0 \iff L_{S_\gvp(A,q)} v = 0, $$
where 
\begin{equation*}
    S_\gvp(A,q) = \bpa{A + 2\nabla \gvp, \, q - (\devp{\gvp}{t} - \gD \gvp  + \abs{\nabla \gvp}^2 + A \cdot \nabla \gvp )} \in L^\infty(Q)^n \times L^\infty(Q) .
\end{equation*}
Thus, for every $(A,q) \in L^\infty(Q)^n \times L^\infty(Q)$, we have
\begin{equation}
    \label{Eq:InvGau}
    \forall \gvp \in G,\quad \gL_{(A,q)} = \gL_{S_\gvp(A,q)}.
\end{equation}
Since $S_\gvp(A,q) \neq (A,q)$ when $ \nabla \gvp \neq 0$, the partial DN map $\gL_{(A,q)}$ is invariant by the gauge transformation of $(A,q)$ and the best we can expect is to prove the recovery of $(A,q)$ up to a gauge transformation $S_\gvp$ for some $\gvp \in G$.
We call the gauge class of $(A,q)$, the equivalence class of $(A,q)$ induced by the equivalence relation defined by
$$ (A,q) \sim (A',q') \iff \exists \gvp \in \gvpSpace,\quad (A',q') = S_\gvp(A,q).$$
From what precedes, we know that the partial DN map is the same for two elements in the same gauge class.

\section{Constructions of geometric optics solutions:} \label{GOS}

Let us consider for this section, $A \in \ASpace$ and $q \in L^\infty(Q)$ along with $\cL_{+} := \devp{}{t} - \gD + A \cdot \nabla + q$ and its formal adjoint $\cL_- := - \devp{}{t} - \gD - A \cdot \nabla + (q - \dive(A))$. 
Our goal in this section is to construct special solutions in $H^{1,2}(Q)$ of the two following problems:
\begin{equation}
    \label{Go1}
        \begin{dcases*}
            \cL_+u = 0 \text{ in } Q, \\
            {u}_{|_{t=0}} = 0 \text{ on } \gO, \\
        \end{dcases*}
        \quad
        \text{ and }\quad 
        \begin{dcases*}
            \cL_- v = 0 \text{ in } Q, \\
            {v}_{|_{t=T}} = 0 \text{ on } \gO. \\
        \end{dcases*}
\end{equation}
More precisely, we will build solutions of \eqref{Go1} of the form 
\begin{equation}
    \label{Go2}
    \forall (t,x) \in Q,\quad \forall \tau > 0,\quad 
    \begin{dcases*}
        u_\tau(t,x) = e^{\gt_\tau(t,x)}\lrbr{ c_+(t,x) + R_{+,\tau}(t,x) }, \\
        v_\tau(t,x) = e^{-\gt_\tau(t,x)}\lrbr{ c_-(t,x) + R_{-,\tau}(t,x) },
    \end{dcases*}
\end{equation}
where the weight function $\gt_\tau$, the principal parts $c_+$ and $c_-$, and, the remainder terms $R_{+,\tau}$ and $R_{-,\tau}$ will be properly defined below.
In the rest of this article, solutions of the form \eqref{Go2} will be called geometric optics solutions.
\new

First, we define the weight function $\gt_\tau$.
For a sufficiently smooth function $\gt_\tau$, one can check that
$$
\begin{dcases*}
    \cL_+ e^{\gt_\tau} = e^{\gt_\tau} \lrbr{ \cL_+ - 2\nabla \gt_\tau \cdot \nabla + (A\cdot \nabla \gt_\tau - \gD \gt_\tau - \abs{\nabla \gt_\tau}^2 + \devp{\gt_\tau}{t})}, \\
    \cL_- e^{-\gt_\tau} = e^{-\gt_\tau} \lrbr{ \cL_- + 2\nabla \gt_\tau \cdot \nabla + (A\cdot \nabla \gt_\tau + \gD \gt_\tau - \abs{\nabla \gt_\tau}^2 + \devp{\gt_\tau}{t})}.
\end{dcases*}
$$
We choose $\gt_\tau$ as a solution of the equation
$$  \devp{\gt_\tau}{t} - \abs{\nabla \gt_\tau}^2 = 0 . $$
Similarly to \cite{choulliLogarithmicStabilityDetermining2018, caroDeterminationConvectionTerms2019, feizmohammadiInverseBoundaryValue2023, kianDeterminingBreakingGauge2024}, we consider weight function of the form 
\begin{equation}
    \label{Go16}
    \gt_\tau(t,x) = \tau^2 t + \tau \psi(x),\quad (t,x) \in \ade{Q},
\end{equation}
where $\psi \in \regu{1}{\ade{\gO}} \cap \regub{2}{\gO}$ and $\tau >0$ denotes a parameter that will be chosen large enough.
With such weights functions, we now have
$$ \cL_\pm e^{\pm\gt_\tau} = e^{\pm\gt_\tau} \cP_{\pm,\tau}, $$
where $\cP_{\pm,\tau} = \cL_\pm + \tau \cJ_\pm + \tau^2 \cI$ with
\begin{equation}
    \label{Go18}
    \begin{dcases*}
        \cI = 1 - \abs{\nabla \psi}^2, \\
        \cJ_+ = - 2 \nabla \psi \cdot \nabla + \lrbr{ A \cdot \nabla \psi - \gD \psi }, \\
        \cJ_- = 2 \nabla \psi \cdot \nabla + \lrbr{ A \cdot \nabla \psi + \gD \psi }.
    \end{dcases*}
\end{equation}
We assume that $\psi$ solves the following eikonal equation:
\begin{equation}
    \label{Hyp2}
    \abs{\nabla \psi}^2 = 1, \quad \text{ on } \gO.
\end{equation}
It will be proved in Appendix \ref{Carleman} that weights functions $\gt_\tau$ of the form \eqref{Go16} with $\psi$ satisfying \eqref{Hyp2}
lead to a Carleman estimate for the operator $\devp{}{t} - \gD + A \cdot \nabla + q$ where $A\in L^\infty(Q)^n$ and $q \in L^\infty(Q)$, when $\psi$ satisfies also
\begin{equation}
    \label{Hyp3}
    \dive(\gD\psi\nabla \psi) \in L^\infty(\gO).
\end{equation} 
We will prove at the same time, a Carleman estimate for the weight $-\gt_\tau$ for the adjoint parabolic operator $-\devp{}{t} - \gD - A \cdot \nabla + q - \dive(A)$ when $\dive(A) \in L^\infty(Q)$.

In this article, we choose to consider the nonlinear weight
\begin{equation}
    \label{Go17}
    \psi :\RR^n \ni x \longmapsto \abs{x-y} \in \RR_+,
\end{equation}
where $y \in \RR^n\backslash\ade{\gO}$ in order to have $\psi$ being smooth in a neighborhood of $\ade{\gO}$ and satisfying \eqref{Hyp2} and \eqref{Hyp3}. 
From now on, we consider $\psi$ satisfying \eqref{Hyp2}.
\new

\subsection{The remainder terms:}

Let us fix $c_{\pm} \in H^{1,2}(Q)$.
The next step is to construct the remainder terms with suitable decay with respect to the large parameter $\tau$.
We get from \eqref{Go1} and \eqref{Go2} that $R_{\pm,\tau}$ should solve in $H^{1,2}(Q)$ the equation
\begin{equation}
    \label{Go3}
    \cP_{\pm,\tau} R_{\pm,\tau} = - \cP_{\pm,\tau} c_\pm.
\end{equation}
Since the problems \eqref{Go1} only involve homogeneous initial or final condition, we can consider the principal parts verifying ${c_+}_{|_{t=0}} = 0$, ${c_-}_{|_{t=T}} = 0$, and define the remainder terms as being solutions of the following problems
\begin{equation}
    \label{Go4}
    \begin{dcases}
        \cP_{+,\tau} R_{+,\tau} = F_+ &\quad \text{ in } Q, \\
        {R_{+,\tau}} = 0 &\quad \text{ on } \gS, \\
        {R_{+,\tau}}(0,\cdot) = 0 &\quad \text{ on } \gO,
    \end{dcases}
\end{equation}
and 
\begin{equation}
    \begin{dcases}
        \label{Go5}
        \cP_{-,\tau} R_{-,\tau} = F_- &\quad \text{ in } Q, \\
        {R_{-,\tau}} = 0 &\quad \text{ on } \gS, \\
        {R_{-,\tau}}(T,\cdot) = 0 &\quad \text{ on } \gO, 
    \end{dcases}
\end{equation}
where $F_\pm = - \cP_{\pm,\tau} c_\pm $.
We will build $R_{\pm,\tau}$ with proper decay property as $\tau \to +\infty$.
More precisely, we prove the following result:
\begin{Proposition}
    \label{Pro:Go1}
    For every $F_\pm \in L^2(Q)$, each of the problem \eqref{Go4} and \eqref{Go5} admits a unique solution in $H^{1,2}(Q)$.
    Moreover, there exist $\tau_0 > 0$ and $C >0$ depending only on $\gO$, $T$, $\psi$, $A$ and $q$ such that for every $\tau > \tau _0$, we get
    \begin{equation}
        \label{Go6}
        \tau \norme{R_{\pm,\tau}}_{L^2(Q)} + \tau^\frac12 \norme{\nabla R_{\pm,\tau}}_{L^2(Q)^n} \leq C \norme{F_{\pm}}_{L^2(Q)}.
    \end{equation}
\end{Proposition}

For this purpose, we use Proposition~\ref{Pro:WP1} and Corollary~\ref{Cor:WP2} with another intermediate lemma which will be stated below. 
We set $q_0 \in L^\infty(Q)$, $A_0 \in L^\infty(Q)^n$, $ L_\pm := \pm \devp{}{t} - \gD \pm A_0 \cdot \nabla + q_0$ along with $P_{\pm,\tau} := L_\pm + \tau J_\pm$ for $\tau > 0$, where $J_\pm := \mp 2 \nabla \psi \cdot \nabla + \lrbr{ A_0 \cdot \nabla \psi \mp \gD \psi }$, such that 
$$ L_\pm e^{\pm\gt_\tau} = e^{\pm\gt_\tau} P_{\pm,\tau} .$$
\begin{Lemme}
    \label{Lem:Go1}
    For every $v_+$, $v_- \in H^{1,2}(Q)$ such as ${v_\pm}_{|_\gS} = 0$, ${v_+}_{|_{t=0}} = 0$ and ${v_-}_{|_{t=T}} = 0$, there exist $\tau_0 > 0$ and $C >0$ which depend only on $\gO$, $T$, $\psi$, $\norme{q_0}_{L^\infty(Q)}$ and $\norme{A_0}_{L^\infty(Q)^n}$ such as for every $\tau > \tau_0$,
    \begin{equation}
        \label{EqPro3}
        \tau \norme{v_\pm}_{L^2(Q)} + \tau^{\frac12} \norme{\nabla v_\pm}_{L^2(Q)^n} \leq C \norme{P_{\pm,\tau} v_\pm }_{L^2(Q)}.
    \end{equation}
\end{Lemme}

\begin{Preuve}
    This result was proved in \cite[Proposition 4.1]{feizmohammadiInverseBoundaryValue2023} when $A_0 = 0$. 
    We extend this analysis to more general convection-diffusion-reaction equation.
    We fix $v_+$, $v_- \in H^{1,2}(Q)$ such as ${v_\pm}_{|_\gS} = 0$, ${v_+}_{|_{t=0}} = 0$ and ${v_-}_{|_{t=T}} = 0$.
    We also consider a constant $\gl > 0$ which will be fixed large enough.
    We set $F_{\pm,\tau} = P_{\pm,\tau} v_\pm$ and we consider lower bound for the quantity
    \begin{equation*}
        \integdd{F_{\pm,\tau} e^{\pm \gl \psi} v_\pm}{Q}{x}{t} = I_\pm + II_\pm + III_\pm + IV_\pm + V_\pm,
    \end{equation*}
    where
    \begin{align*}
        I_\pm &= \integdd{q_0v_\pm^2e^{\pm\gl\psi}}{Q}{x}{t}, \quad II_\pm = -\integdd{\lrpa{\gD v_\pm} v_\pm e^{\pm\gl\psi} }{Q}{x}{t}, \\
        III_\pm &= \pm\integdd{\lrpa{A_0\cdot \nabla v_\pm }v_\pm e^{\pm\gl\psi} }{Q}{x}{t},\quad IV_\pm = \pm \integdd{\lrpa{\devp{v_\pm}{t}} v_\pm e^{\pm \gl\psi}}{Q}{x}{t}, \\
        V_\pm &= \tau\integdd{\lrpa{J_\pm v_\pm} v_\pm e^{\pm\gl\psi}}{Q}{x}{t}.
    \end{align*}
    In all the remaining parts of the proof, we set $T_+ = T$ and $T_- = 0$.
    Note first that
    \begin{equation}
        \label{IEF1}
        I_\pm \geq -\norme{q_0}_{L^\infty(Q)}\integdd{v_\pm^2 e^{\pm\gl\psi}}{Q}{x}{t}.
    \end{equation}
    Using the fact that ${v_\pm}_{|_\gS} = 0$, we integrate by parts to obtain,
    \begin{align*}
        II_\pm &= \integdd{\nabla v_\pm \cdot \nabla\lrpa{e^{\pm\gl\psi} v_\pm }}{Q}{x}{t} \\
        &= \integdd{\abs{\nabla v_\pm}^2 e^{\pm\gl\psi}}{Q}{x}{t} + \frac12 \integdd{\nabla v_\pm^2 \cdot \nabla e^{\pm\gl\psi}}{Q}{x}{t} \\
        &= \integdd{\abs{\nabla v_\pm}^2 e^{\pm\gl\psi}}{Q}{x}{t} - \frac12 \integdd{v_\pm^2 \gD e^{\pm\gl\psi}}{Q}{x}{t} \\
        &= \integdd{\abs{\nabla v_\pm}^2 e^{\pm\gl\psi}}{Q}{x}{t} - \frac12 \integdd{v_\pm^2 e^{\pm\gl\psi}\lrbr{\gl^2 \pm \gl \gD \psi} }{Q}{x}{t},
    \end{align*}
    from which we derive
    \begin{equation}
        \label{IEF2}
        II_\pm \geq \integdd{\abs{\nabla v_\pm}^2 e^{\pm\gl\psi}}{Q}{x}{t} - (\gl^2 + \gl \norme{\gD\psi}_{L^\infty(Q)}) \integdd{v_\pm^2 e^{\pm\gl\psi}}{Q}{x}{t} .
    \end{equation}
    We consider $h > 0$ that will be fixed later and we apply Young inequality to get
    \begin{align*}
        III_\pm &\geq -\frac12 \lrpa{ h\integdd{v_\pm^2 e^{\pm\gl\psi}}{Q}{x}{t} + \frac{1}{h}\integdd{\lrpa{A_0\cdot \nabla v_\pm}^2 e^{\pm\gl\psi}}{Q}{x}{t} } \\
        &\geq -\frac12\lrpa{h\integdd{v_\pm^2 e^{\pm\gl\psi}}{Q}{x}{t} + \frac{1}{h}\norme{A_0}^2_{L^\infty(Q)^n} \integdd{\abs{\nabla v_\pm}^2 e^{\pm\gl\psi}}{Q}{x}{t} }.
    \end{align*}
    Fixing $h = \norme{A_0}^2_{L^\infty(Q)^n} +1 $, we obtain
    \begin{equation}
        \label{IEF3}
        III_\pm \geq -\frac12\lrpa{\lrpa{\norme{A_0}^2_{L^\infty(Q)^n} +1}\integdd{v_\pm^2 e^{\pm\gl\psi}}{Q}{x}{t} + \integdd{\abs{\nabla v_\pm}^2 e^{\pm\gl\psi}}{Q}{x}{t} }.
    \end{equation}
    Since ${v_\pm}_{|_{t=T_\pm}} = 0$, we get
    $$ IV_\pm = \pm \frac12 \integdd{\devp{v_\pm^2}{t} e^{\pm\gl\psi}}{Q}{x}{t} = \frac12 \integd{v_\pm^2(T_\mp,x) e^{\pm\gl\psi} }{\gO}{x}, $$
    which implies
    \begin{equation}
        \label{IEF4}
        IV_\pm \geq 0.
    \end{equation}
    Finally, using the fact that ${v_\pm}_{|_\gS} = 0$, we have,
    \begin{align*}
        V_\pm &= \tau \integdd{A_0\cdot \nabla \psi v_\pm^2 e^{\pm\gl\psi} }{Q}{x}{t} \mp \tau \integdd{\lrpa{ 2\nabla \psi \cdot \nabla v_\pm + \gD \psi v_\pm} v_\pm e^{\pm\gl\psi} }{Q}{x}{t} \\
        &\geq -\tau \norme{A_0}_{L^\infty(Q)^n} \integdd{v_\pm^2 e^{\pm\gl\psi}}{Q}{x}{t} \mp \tau \lrbr{\integdd{\lrpa{\nabla \psi \cdot \nabla v_\pm^2} e^{\pm\gl\psi} }{Q}{x}{t} + \integdd{\gD \psi v_\pm^2 e^{\pm\gl\psi} }{Q}{x}{t}} \\
        &\geq -\tau \norme{A_0}_{L^\infty(Q)^n} \integdd{v_\pm^2 e^{\pm\gl\psi}}{Q}{x}{t} \mp \tau \lrbr{ - \integdd{v_\pm^2 \dive\lrpa{e^{\pm\gl\psi} \nabla \psi} }{Q}{x}{t} + \integdd{\gD \psi v_\pm^2 e^{\pm\gl\psi} }{Q}{x}{t}} \\
        &\geq -\tau \norme{A_0}_{L^\infty(Q)^n} \integdd{v_\pm^2 e^{\pm\gl\psi}}{Q}{x}{t} \mp \tau \lrbr{ \mp \gl \integdd{v_\pm^2 e^{\pm\gl\psi}\abs{\nabla \psi}^2 }{Q}{x}{t} }
    \end{align*}
    and we deduce that 
    \begin{equation}
        \label{IEF5}
        V_\pm \geq \tau \gl \integdd{v_\pm^2 e^{\pm\gl\psi}}{Q}{x}{t} -\tau \norme{A_0}_{L^\infty(Q)^n} \integdd{v_\pm^2 e^{\pm\gl\psi}}{Q}{x}{t}.
    \end{equation}
    Combining \eqref{IEF1}-\eqref{IEF5}, for 
    $$ K_{\tau,\gl} := \lrpa{\gl - \norme{A_0}_{L^\infty(Q)^n}} \tau - C_-(1+\gl+\gl^2) ,$$
    with $C_- > 0$ which depends only on $\gO$, $T$, $\psi$, $\norme{q_0}_{L^\infty(Q)}$ and $\norme{A_0}_{L^\infty(Q)^n}$, we obtain
    \begin{equation}
        \label{IEF6}
        \integdd{F_{\pm,\tau} e^{\pm \gl \psi} v_\pm}{Q}{x}{t} \geq \frac{1}{2}\integdd{\abs{\nabla v_\pm}^2 e^{\pm\gl\psi}}{Q}{x}{t} + K_{\tau,\gl} \integdd{v_\pm^2 e^{\pm\gl\psi}}{Q}{x}{t}.
    \end{equation}
    We now give a suitable lower bound of $K_{\tau,\gl}$ in order to conclude this proof. 
    We set $\ga := 7C_-$ and $\gl_0 := \max\ense{2;4\norme{A_0}_{L^\infty(Q)^n}}$ and notice that for every $\gl \geq \gl_0$, we have
    $$ \frac34C_-\gl^2 - C_-\gl - C_- \geq 0 . $$
    Since $\frac34 C_- = \frac\ga 4 - C_-$, for every $\gl \geq \gl_0$ and $\tau \geq \ga \gl$, we find
    \begin{align*}
        C_-\lrpa{1+\gl+\gl^2} &\leq \frac{\gl^2}4\ga \leq \lrpa{\frac\gl2 - \norme{A_0}_{L^\infty(Q)^n}} \ga \gl \\
        &\leq \lrpa{\frac\gl2 - \norme{A_0}_{L^\infty(Q)^n}}\tau.
    \end{align*}
    It follows that, for every $\gl \geq \gl_0$ and $\tau \geq \ga \gl$, we have $\tau\frac{\gl}{2} \leq K_{\tau,\gl}$ and with \eqref{IEF6}, we get
    \begin{equation}
        \label{IEF7}
        \integdd{F_{\pm,\tau} e^{\pm \gl \psi} v_\pm}{Q}{x}{t} \geq \frac{1}{2}\integdd{\abs{\nabla v_\pm}^2 e^{\pm\gl\psi}}{Q}{x}{t} + \tau\frac{\gl}{2} \integdd{v_\pm^2 e^{\pm\gl\psi}}{Q}{x}{t}.
    \end{equation}
    Now using Young inequality, we have 
    $$  \integdd{F_{\pm,\tau} e^{\pm \gl \psi} v_\pm}{Q}{x}{t}  \leq \frac12 \lrpa{ \frac{1}{\tau}  \integdd{F_{\pm,\tau}^2 e^{\pm \gl \psi}}{Q}{x}{t} + \tau  \integdd{v_\pm^2 e^{\pm \gl \psi}}{Q}{x}{t} },
    $$
    for every $\tau >0$ and applying \eqref{IEF7}, for every $\gl \geq \gl_0$ and $\tau \geq \ga \gl$, we get
    \begin{equation*}
        \frac{1}{2\tau}\integdd{F_{\pm,\tau}^2 e^{\pm \gl \psi}}{Q}{x}{t} \geq \frac{1}{2}\integdd{\abs{\nabla v_\pm}^2 e^{\pm\gl\psi}}{Q}{x}{t} + \frac{\tau}{2} \lrpa{\gl - 1}\integdd{v_\pm^2 e^{\pm\gl\psi}}{Q}{x}{t}.
    \end{equation*}
    We note $M := \max_\ade{\gO} \abs{\psi}$.
    For every $\gl \geq \gl_0$, we have $e^{-\gl M} \leq e^{\pm \gl \psi} \leq e^{\gl M}$, and for $\tau \geq \ga \gl$, we obtain 
    \begin{equation*}
        e^{2\gl M}\norme{F_{\pm,\tau}}_{L^2(Q)}^2 \geq \tau \norme{\nabla v_\pm}_{L^2(Q)^n}^2 + \tau^2\lrpa{\gl - 1}\norme{v_\pm}_{L^2(Q)}^2.
    \end{equation*}
    In order to conclude, we arbitrarily fix $\gl \geq \gl_0$, define $\tau_0 := \ga \gl$ and deduce the existence of $C$ positive which depends only on $\gO$, $T$, $\psi$, $\norme{q_0}_{L^\infty(Q)}$ and $\norme{A_0}_{L^\infty(Q)^n}$ such that we have \eqref{EqPro3} for every $\tau > \tau_0$. 
    We finally notice that by construction, $\tau_0$ depends only on $\gO$, $T$, $\psi$, $\norme{q_0}_{L^\infty(Q)}$ and $\norme{A_0}_{L^\infty(Q)^n}$.
\end{Preuve}

\begin{Preuve}[of Proposition~\ref{Pro:Go1}]
    Let us fix $\tau > 0$.
    Since 
    $$ \cP_{\pm,\tau} = \cL_\pm + \tau \lrpa{\mp 2 \nabla \psi \cdot \nabla + \lrbr{ A \cdot \nabla \psi \mp \gD \psi }} , $$
    we deduce from Proposition~\ref{Pro:WP1} and Corollary~\ref{Cor:WP2} that for every $F_\pm \in L^2(Q)$, each of the problem \eqref{Go4} and \eqref{Go5} admits a unique solution $R_{\pm,\tau}$ in $H^{1,2}(Q)$.
    \new

    Now for $F_{\pm}$ fixed in $L^2(Q)$, using Lemma~\ref{Lem:Go1} with $v_\pm = R_{\pm,\tau}$ and $L_\pm = \cL_\pm$, we get for every $\tau > 0$, that $\cP_{\pm,\tau} = P_{\pm,\tau}$ and that there exist $\tau_0 > 0$ and $C >0$ which depend only on $\gO$, $T$, $\psi$, $\norme{q}_{L^\infty(Q)}$, $A$ such as for every $\tau > \tau_0$, we have
    \begin{equation*}
        \tau \norme{R_{\pm,\tau}}_{L^2(Q)} + \tau^{\frac12} \norme{\nabla R_{\pm,\tau}}_{L^2(Q)^n} \leq C \norme{\cP_{\pm,\tau} R_{\pm,\tau}}_{L^2(Q)} = C \norme{F_\pm}_{L^2(Q)}.
    \end{equation*}
    This completes the proof of the proposition.
\end{Preuve}

\subsection{The principal parts:}

We complete this section by constructing the principal parts $c_\pm$ of the geometric optics solutions \eqref{Go2}. 
Using Proposition~\ref{Pro:Go1} and \eqref{Go3}, we know that we can obtain a convenient decay property for the remainder terms if 
$$ \cP_{\pm,\tau} c_\pm = \cL_{\pm} c_\pm \quad \text{ in } Q.$$
To get the latest equality, knowing that $\cP_{\pm,\tau} = \cL_\pm + \tau\cJ_\pm$ with $\cJ_\pm$ as in \eqref{Go18}, we consider principal parts $c_+$ and $c_-$ solving the following transports equations 
\begin{equation}
    \label{Go7}
    \begin{dcases}
        \cJ_+ c_+ = 0 &\text{ in } Q, \\
        c_+(0,\cdot) = 0 &\text{ on } \gO,
    \end{dcases}
\end{equation} 
and 
\begin{equation}
    \label{Go8}
    \begin{dcases}
        \cJ_- c_- = 0 &\text{ in } Q, \\
        c_-(T,\cdot) = 0 &\text{ on } \gO.
    \end{dcases}
\end{equation} 
For the rest of this section, we consider $\psi$ as in \eqref{Go17} where $y \in \RR^n\backslash\ade{\gO}$.
From now on, we denote by $\SS^{n-1}$ the unit sphere in $\RR^n$ and we state the following result.

\begin{Lemme}
    \label{Lem:Go2}
    Let $y \in B(x_0,\gve)$ fixed. 
    If $\tilde{A} \in \regub{1}{0,T;\regu{0}{\RR^n}} \cap \regub{0}{0,T;\regu{2}{\RR^n}}$, $h \in \reguc{\infty}{\SS^{n-1}}$ and $\chi \in \reguc{\infty}{0,T}$ then the map $c$ defined for every $(t,x) \in \inteoo{0}{T} \times \RR^n\backslash\ense{y}$ by
    \begin{equation}
        \label{Go9}
        c(t,x) = \exp\lrpa{ \Integd{\frac{\tilde{A}(t,y+s(x-y))}{2} \cdot (x-y)}{0}{1}{s}}\abs{x-y}^{-\frac{n-1}{2}}h\lrpa{\frac{x-y}{\abs{x-y}}}\chi(t),
    \end{equation}
    verifies $- 2 \nabla \psi \cdot \nabla c + \lrbr{\tilde{A} \cdot \nabla \psi - \gD \psi} c = 0$ over $Q$ with $c(0,\cdot) = c(T,\cdot) = 0$ on $\gO$ and $c_{|_Q} \in H^{1,2}(Q)$.
\end{Lemme}

For the construction of the principal parts $c_\pm$ of the geometric optics solutions \eqref{Go2} using Lemma~\ref{Lem:Go2}, we need first to justify the existence of an extension of $A$ in $\regub{1}{0,T;\regu{0}{\RR^n}} \cap \regub{0}{0,T;\regu{2}{\RR^n}}$.
For this purpose, we consider a bounded extension operator $P : \regu{0}{\ade{\gO}} \rightarrow \regu{0}{\RR^n}$ such that $P$ induces a bounded operator from $\regu{2}{\ade{\gO}}$ to $\regu{2}{\RR^n}$. 
Such operator exists thanks to the Seeley extension operator \cite{seeleyExtensionFunctionsDefined1964} and the fact that $\gO$ is a $\cC^2$ open set.
We then set
$$ \forall (t,x) \in \inteoo{0}{T}\times \RR^n ,\quad \tilde{A}(t,x) = P(A(t,\cdot))(x) ,$$
in order to conclude the existence. 
We now prove Lemma \ref{Lem:Go2}.

\begin{Preuve}[of Lemma~\ref{Lem:Go2}]
    To prove this result, we will construct solutions $\mu$ of 
    \begin{equation}
        \label{Go10}
        \begin{dcases}
            \cJ \mu = 0 &\text{ in } \inteoo{0}{T} \times \RR^n \backslash \ense{y}, \\
            \mu(0,\cdot) = 0 &\text{ on } \RR^n \backslash \ense{y}, \\
            \mu(T,\cdot) = 0 &\text{ on } \RR^n \backslash \ense{y},
        \end{dcases}
    \end{equation} 
    where the operator $\cJ$ is defined by $\cJ = - 2 \nabla \psi \cdot \nabla + \lrbr{\tilde{A} \cdot \nabla \psi - \gD \psi}$ and $\tilde{A}$ is a continuous function in $\inteoo{0}{T} \times\RR^n$. 
    For this purpose, we set our transport equation into polar coordinates centered at $y$ and consider the following maps
    \begin{equation*}
        \begin{dcases*}
            \phi : \RR_+^* \times \SS^{n-1} \ni (r,\gt) \longmapsto y+ r\gt \in \RR\backslash\ense{y}, \\
            \Phi : \inteoo{0}{T} \times \RR_+^* \times \SS^{n-1} \ni (t,r,\gt) \longmapsto (t,\phi(r,\gt)) \in \inteoo{0}{T} \times \RR\backslash\ense{y}.
        \end{dcases*}
    \end{equation*}
    Now, for a map $u$ defined on $\RR^n$ (resp. $\inteoo{0}{T} \times \RR^n$), we write $\gk u = u \circ \phi$ (resp. $Ku = u \circ \Phi$).
    We first notice that for every $(r,\gt) \in \RR^*_+ \times \SS^{n-1}$,
    \begin{equation}
        \label{Go11}
        \gk\lrpa{\nabla \psi}(r,\gt) = \gt \quad\text{ and }\quad \gk(\gD\psi)(r,\gt) = \frac{n-1}{r}.
    \end{equation}
    Moreover, for $u \in \regu{1}{\RR^n}$, we can compute that
    \begin{equation}
        \label{Go12}
        \forall (r,\gt) \in \RR^*_+ \times \SS^{n-1},\quad \kappa(2\nabla \psi \cdot \nabla u)(r,\gt) = 2\devp{(\gk u)}{r}(r,\gt) .
    \end{equation}
    
    Combining \eqref{Go11} and \eqref{Go12}, solving problem \eqref{Go10} becomes equivalent to
    \begin{equation}
        \label{Go13}
        \begin{dcases}
            \devp{\tilde{\mu}}{r}(t,r,\gt) = \frac{1}{2}\lrbr{ K\tilde{A}(t,r,\gt) \cdot \gt - \frac{n-1}{r} }\tilde{\mu}(t,r,\gt) &\text{ for } (t,r,\gt) \in \inteoo{0}{T} \times \RR_+^*\times \SS^{n-1}, \\
            \tilde{\mu}(0,\cdot) = 0 &\text{ for } (r,\gt) \in \RR_+^*\times \SS^{n-1}, \\
            \tilde{\mu}(T,\cdot) = 0 &\text{ for } (r,\gt) \in \RR_+^*\times \SS^{n-1},
        \end{dcases}
    \end{equation} 
    when solutions are considered regular enough. 
    We then deduce that for every $h \in \reguc{\infty}{\SS^{n-1}}$ and $\chi \in \reguc{\infty}{0,T}$, the map $\tilde{\mu}$ defined on $\inteoo{0}{T}\times \RR_+^* \times \SS^{n-1}$ by 
    \begin{equation}
        \label{Go14}
        \tilde{\mu}(t,r,\gt) = \exp\lrpa{ \Integd{\frac{K\tilde{A}(t,s,\gt)}{2} \cdot \gt}{0}{r}{s}}r^{-\frac{n-1}{2}}h(\gt)\chi(t),
    \end{equation}
    is a solution of \eqref{Go13}.
    Moreover, assuming that $\tilde{A} \in \regub{1}{0,T;\regu{0}{\RR^n}} \cap \regub{0}{0,T;\regu{2}{\RR^n}}$, we get $\lrpa{\tilde{\mu} \circ \Phi^{-1}}_{|_Q} \in H^{1,2}(Q)$ with \eqref{Hyp1}. 
    The change of variable $s = ur $ in \eqref{Go14} allows us to conclude that $c$ defined as in \eqref{Go9} is a solution of \eqref{Go10} over $\inteoo{0}{T}\times \RR^n\backslash\ense{y}$ which satisfies $c_{|_Q} \in H^{1,2}(Q)$.
\end{Preuve}

\section{Proof of Theorem \ref{Theo1}:}			\label{UnPr}

For the rest of this section, we fix $(A_1,q_1)$ and $(A_2,q_2)$ two elements of $\cA \times L^\infty(Q)$ such that ${A_1}_{|_\gS} = {A_2}_{|_\gS}$ and assume that $\gL_{(A_1,q_1)} = \gL_{(A_2,q_2)}$. 
We divide the proof of Theorem~\ref{Theo1} in different steps.
In the first step, we introduce direct consequences of the assumption mentioned above. 
We then use results from Section \ref{GOS} to obtain an integral identity related to a parameter $\tau$, the principal parts and the remainder terms of geometric optics solutions. 
We also use a Carleman estimate to derive asymptotic properties of the integral identity when $\tau$ goes to infinity.
The second step is devoted to results derived from the integral identity. 
More specifically, we first prove that we can determine $q$ from the knowledge of $\gL_{(A,q)}$ and $A$ (Theorem~\ref{Theo2}), and finally, we prove that the X-ray transform of the continuous zero-extension of $A_2-A_1$ vanishes over every line that intersects $B(x_0,\gve)$, after recalling the definition of the X-ray transform of a compactly supported vector field.
In the last step, we construct $\gvp$ and conclude using the gauge invariance.

\subsection{Notations and first results:}
\label{Nota}

For the convenience of readers, we set $A_+ := A_1$ and $A_- := A_2$.
For $j=1,2$, let us consider $w_j$ in $\SolSpace$ such as 
\begin{equation}
    \label{Mt1}
        \begin{dcases}
            \devp{w_j}{t} - \gD w_j + A_j \cdot \nabla w_j + q_jw_j = 0 &\text{ in } Q, \\
            {w_j}_{|_{t=0}} = 0 &\text{ on } \gO,
        \end{dcases}
\end{equation}
and $w_1 = w_2$ on $\gS$. 
Let us also consider $v \in \SolSpace$ such as 
\begin{equation}
    \label{Mt2}
        \begin{dcases}
            - \devp{v}{t} - \gD v - A_2 \cdot \nabla v + (q_2 - \dive(A_2))v = 0 &\text{ in } Q, \\
            {v}_{|_{t=T}} = 0 &\text{ on } \gO.
        \end{dcases}
\end{equation}
Let $w := w_1-w_2$ and observe that $w$ satisfies the following properties,
\begin{equation}
    \label{Mt3}
        \begin{dcases}
            \devp{w}{t} - \gD w + A_2 \cdot \nabla w + q_2 w = (A_2-A_1)\cdot \nabla w_1 + (q_2-q_1) w_1 &\text{ in } Q, \\
            w_{|_\gS} = 0 &\text{ on } \gS \\
            \devp{w}{\nu}_{|_{\gS_-\lrpa{x_0,2\gve}}} = 0 &\text{ on } \gS, \\
            {w}_{|_{t=0}} = 0 &\text{ on } \gO.
        \end{dcases}
\end{equation}
Multiplying the equation \eqref{Mt3} by $v$, integrating by parts and using \eqref{Mt2}, we obtain 
\begin{align*}
    \integdd{\lrbr{(A_2-A_1)\cdot \nabla w_1 + (q_2-q_1) w_1}v}{Q}{x}{t} =& \integdd{\lrbr{\devp{w}{t} - \gD w + A_2 \cdot \nabla w + q_2 w}v}{Q}{x}{t} \\
    =& \integdd{\lrbr{ -\devp{v}{t} -\gD v - A_2 \cdot \nabla v + (q_2-\dive(A_2)) v}w}{Q}{x}{t}  \\
    &+ \integdd{\lrpa{w\devp{v}{\nu} - v\devp{w}{\nu} + v w A_2\cdot \nu} }{\gS}{\gs(x)}{t} \\
    =& - \integdd{v\devp{w}{\nu}}{\gS}{\gs(x)}{t} .
\end{align*}
Therefore, we get
\begin{equation}
    \label{Mt4}
    \integdd{\lrbr{(A_2-A_1)\cdot \nabla w_1 + (q_2-q_1) w_1}v}{Q}{x}{t} = - \integdd{v\devp{w}{\nu}}{\gS\backslash\gS_-(x_0,2\gve)}{\gs(x)}{t} .
\end{equation}
\new

Let $\cL_{+} := \devp{}{t} - \gD + A_1 \cdot \nabla + q_1$ and $\cL_- := - \devp{}{t} - \gD - A_2 \cdot \nabla + (q_2 - \dive(A_2))$. 
We fix $y \in B(x_0,\gve)$ and define the map $\psi: \RR^n \ni x \longmapsto \abs{x-y}\in \RR_+$ along with $\gt_\tau$ as in \eqref{Go16} for every $\tau > 0$.
For such $y$ and $\tau$, we also define the operators 
\begin{equation}
    \label{Mt5}
    \begin{dcases*}
        \cJ_+ = - 2 \nabla \psi \cdot \nabla + \lrbr{ A_1 \cdot \nabla \psi - \gD \psi }, \\
        \cJ_- = 2 \nabla \psi \cdot \nabla + \lrbr{ A_2 \cdot \nabla \psi + \gD \psi },
    \end{dcases*}
\end{equation}
and the operator $\cP_{\pm,\tau} := \cL_\pm + \tau \cJ_\pm$, which verifies 
$$ \cL_\pm e^{\pm\gt_\tau} = e^{\pm\gt_\tau} \cP_{\pm,\tau} .$$

We consider for $i=1,2$, $\tilde{A}_i$ an extension of $A_i$ in $\regub{1}{0,T;\regu{0}{\RR^n}} \cap \regub{0}{0,T;\regu{2}{\RR^n}}$.
For the convenience of readers, we set $\tilde{A}_+ := \tilde{A}_1$ and $\tilde{A}_- := \tilde{A}_2$. 
For a fixed $h_\pm \in \reguc{\infty}{\SS^{n-1}}$ and $\chi_\pm \in \reguc{\infty}{0,T}$ and for every $(t,x) \in \inteoo{0}{T} \times \RR^n\setminus\ense{y}$, we define
\begin{equation}
    \label{Mt6}
    c_\pm(t,x) = \exp\lrpa{ \pm \Integd{\frac{\tilde{A}_\pm(t,y+s(x-y))}{2} \cdot (x-y)}{0}{1}{s}}\abs{x-y}^{-\frac{n-1}{2}}h_\pm\lrpa{\frac{x-y}{\abs{x-y}}}\chi_\pm(t).
\end{equation}
Notice that $c_\pm$ depend only on $y$, $\tilde{A}_\pm$, $h_\pm$, and $\chi_\pm$.
In view of Lemma~\ref{Lem:Go2}, ${c_\pm}_{|_Q}$ belong to $\SolSpace$ and fulfill the conditions
\begin{equation}
    \label{Mt7}
    \begin{dcases}
        \cJ_\pm {c_\pm} = 0 &\text{ in } Q, \\
        {c_\pm}(0,\cdot) = 0 &\text{ on } \gO, \\
        {c_\pm}(T,\cdot) = 0 &\text{ on } \gO.
    \end{dcases}
\end{equation}
We introduce also,
\begin{equation}
    \label{Mt8}
    \forall (t,x) \in Q,\quad \forall \tau > 0,\quad 
    \begin{dcases*}
        u_\tau(t,x) = e^{\tau^2t+\tau \psi(x)}\lrbr{ c_+(t,x) + R_{+,\tau}(t,x) }, \\
        v_\tau(t,x) = e^{-\tau^2t-\tau\psi(x)}\lrbr{ c_-(t,x) + R_{-,\tau}(t,x) },
    \end{dcases*}
\end{equation}
where, according to Proposition~\ref{Pro:Go1}, $R_{+,\tau}$ and $R_{-,\tau}$ are defined as the unique element in $\SolSpace$ verifying,
\begin{equation}
    \begin{dcases}
        \cP_{\pm,\tau} R_{\pm,\tau} = - \cP_{\pm,\tau} c_\pm &\quad \text{ in } Q, \\
        {R_{\pm,\tau}} = 0 &\quad \text{ on } \gS, \\
        {R_{+,\tau}}(0,\cdot) = 0 &\quad \text{ on } \gO, \\
        {R_{-,\tau}}(T,\cdot) = 0 &\quad \text{ on } \gO.
    \end{dcases}
\end{equation}
By definition, we know that $u_\tau$ and $v_\tau$ are in $\SolSpace$ with
\begin{equation}
    \label{Mt9}
        \begin{dcases*}
            \cL_+ u_\tau = 0 \quad \text{ in } Q, \\
            {u_\tau}_{|_{t=0}} = 0 \quad \text{ on } \gO, \\
        \end{dcases*}
        \quad \text{ and }\quad 
        \begin{dcases*}
            \cL_- v_\tau = 0 \quad\text{ in } Q, \\
            {v_\tau}_{|_{t=T}} = 0 \quad\text{ on } \gO, \\
        \end{dcases*}
\end{equation}
and $u_\tau$ (resp. $v_\tau$) depends only on $y$, $\tau$, $\tilde{A}_+$, $h_+$ and $\chi_+$ (resp. $\tilde{A}_-$, $h_-$ and $\chi_-$).
We also know from Proposition~\ref{Pro:Go1} that there exist $\tau_0 > 0$ and $C >0$ depending only on $\gO$, $T$, $A_1$, $A_2$, $q_1$ and $q_2$ such that for every $\tau > \tau _0$, we have
\begin{equation}
    \label{Mt10}
    \tau \norme{R_{\pm,\tau}}_{L^2(Q)} + \tau^\frac12 \norme{\nabla R_{\pm,\tau}}_{L^2(Q)^n} \leq C \norme{\cL_{\pm}c_\pm}_{L^2(Q)},
\end{equation}
since $\cL_\pm c_\pm = \cP_{\pm,\tau} c_\pm$ according to \eqref{Mt5}.
\new

We state the following Carleman estimate, derived from a more general result whose proof is postponed to Appendix \ref{Carleman}.
\begin{Theoreme}
    \label{The:ICA}
    Let $A \in L^\infty(Q)^n$ and $q \in L^\infty(Q)$.
    For every $v \in H^{1,2}(Q)$ satisfying $v_{|_\gS} = 0$ and $v_{|_{t=0}} = 0$, there exist $\tau_1 > 0$ and $C' >0$ which depend only on $\gO$, $T$, $\norme{q}_{L^\infty(Q)}$ and $\norme{A}_{L^\infty(Q)^n}$ such that for every $\tau > \tau_1$, we have
    \begin{multline}
        \label{Mt37}
        \tau \integdd{\abs{e^{-\gt_\tau}\devp{v}{\nu}}^2 \abs{\devp{\psi}{\nu}}}{\gS_+(y,0)}{\gs(x)}{t} + \tau^2 \integdd{\abs{e^{-\gt_\tau}v}^2}{Q}{x}{t} + \tau \integdd{\abs{\nabla(e^{-\gt_\tau}v)}^2}{Q}{x}{t} \\ \leq C'\lrbr{ \integdd{\abs{e^{-\gt_\tau} Lv}^2}{Q}{x}{t} + \tau \integdd{\abs{e^{-\gt_\tau}\devp{v}{\nu}}^2 \abs{\devp{\psi}{\nu}}}{\gS_-(y,0)}{\gs(x)}{t} },
    \end{multline}
    where $L= \devp{}{t} - \gD + A\cdot \nabla + q$.
\end{Theoreme}
The next step is to consider \eqref{Mt4} with the previously built geometric optics solutions along with the Carleman estimate \eqref{Mt37}.
From \eqref{Mt9}, we obtain 
\begin{equation}
    \label{Mt11}
    \integdd{\lrbr{(A_2-A_1)\cdot \nabla u_\tau + (q_2-q_1) u_\tau}v_\tau}{Q}{x}{t} = - \integdd{v_\tau\devp{w}{\nu}}{\gS\backslash\gS_-(x_0,2\gve)}{\gs(x)}{t},
\end{equation}
where $w = u_\tau - w_2$ with $w_2$ in $\SolSpace$ being the unique solution of
\begin{equation*}
    \begin{dcases*}
        \devp{w_2}{t} - \gD w_2 + A_2 \cdot \nabla w_2 + q_2w_2 = 0, \\
        {w_2}_{|_\gS} = {u_\tau}_{|_\gS}, \\
        {w_2}_{|_{t=0}} = 0. \\
    \end{dcases*}
\end{equation*}
We recall from \eqref{Mt3} that $w$ satisfies
\begin{equation}
    \label{Mt12}
    \begin{dcases*}
        \devp{w}{t} - \gD w + A_2 \cdot \nabla w + q_2 w = (A_2-A_1)\cdot \nabla u_\tau + (q_2-q_1) u_\tau, \\
        w_{|_\gS} = 0, \\
        \devp{w}{\nu}_{|_{\gS_-\lrpa{x_0,2\gve}}} = 0, \\
        {w}_{|_{t=0}} = 0. \\
    \end{dcases*}
\end{equation}
\new

We consider the asymptotic properties of \eqref{Mt11} when $\tau$ goes to infinity.
Using the definition of $u_\tau$ and $v_\tau$, we derive 
\begin{equation}
    \integdd{\lrbr{(A_2-A_1)\cdot \nabla u_\tau + (q_2-q_1) u_\tau}v_\tau}{Q}{x}{t} =  \tau I_\tau + II_\tau,
\end{equation}
with 
\begin{equation}
    \begin{dcases}
        \label{Mt13}
        I_\tau = \integdd{\lrbr{(A_2-A_1)\cdot \nabla \psi}(c_+ + R_{+,\tau})(c_- + R_{-,\tau})}{Q}{x}{t}, \\
        II_\tau = \integdd{\lrbr{(A_2-A_1)\cdot \nabla (c_+ + R_{+,\tau}) + (q_2-q_1) (c_+ + R_{+,\tau})}(c_- + R_{-,\tau})}{Q}{x}{t}.
    \end{dcases} 
\end{equation}
Applying the decay property \eqref{Mt10}, we have the following limits 
\begin{equation}
    \label{Mt14}
    \begin{dcases}
        I_\tau \xrightarrow[\tau\to\infty]{} \integdd{\lrbr{(A_2-A_1)\cdot \nabla \psi}c_+ c_- }{Q}{x}{t}, \\
        II_\tau \xrightarrow[\tau\to\infty]{} \integdd{\lrbr{(A_2-A_1)\cdot \nabla c_+ + (q_2-q_1) c_+ }c_- }{Q}{x}{t}.
    \end{dcases} 
\end{equation}

Since $y \in B(x_0,\gve)$, we get
$$ \forall x \in \gG, \quad (x-x_0) \cdot \nu(x) \leq \abs{x_0-y} + (x-y) \cdot \nu(x), $$
which implies $\gS_-(y,\gve) \subset \gS_-(x_0,2\gve)$.
Therefore, we find 
\begin{equation}
    \label{Mt15}
    \abs{\integdd{v_\tau\devp{w}{\nu}}{\gS\backslash\gS_-(x_0,2\gve)}{\gs(x)}{t}} \leq \integdd{\abs{v_\tau\devp{w}{\nu}}}{\gS\backslash\gS_-(y,\gve)}{\gs(x)}{t}.
\end{equation}
Since $R_{-,\tau}$ vanishes on $\gS$, we have 
$$ \integdd{\abs{v_\tau\devp{w}{\nu}}}{\gS\backslash\gS_-(y,\gve)}{\gs(x)}{t} = \integdd{e^{- {\gt_\tau} }\abs{c_-\devp{w}{\nu}}}{\gS\backslash\gS_-(y,\gve)}{\gs(x)}{t}. $$
Applying Cauchy-Schwartz inequality to the right-hand side, we get 
$$ \integdd{\abs{v_\tau\devp{w}{\nu}}}{\gS\backslash\gS_-(y,\gve)}{\gs(x)}{t} \leq \norme{ c_- }_{L^2(\gS)} \lrpa{\integdd{e^{-2{\gt_\tau} }\abs{\devp{w}{\nu}}^2}{\gS\backslash\gS_-(y,\gve)}{\gs(x)}{t}}^{\frac12} . $$
Using the definition of $\gS\backslash\gS_-(y,\gve)$ and the fact that $\devp{\psi}{\nu}(x) = \frac{x-y}{\abs{x-y}}\cdot \nu(x) $ for every $x$ in $\gG$, we obtain for $m:=\inf_\gG \psi^{-1} > 0$,
$$ \integdd{\abs{v_\tau\devp{w}{\nu}}}{\gS\backslash\gS_-(y,\gve)}{\gs(x)}{t} \leq \norme{ c_- }_{L^2(\gS)} (m\gve)^{-\frac12}\lrpa{\integdd{e^{-2{\gt_\tau} }\abs{\devp{w}{\nu}}^2 \abs{\devp{\psi}{\nu}}}{\gS\backslash\gS_-(y,\gve)}{\gs(x)}{t}}^{\frac12} . $$
Having $\gS\backslash\gS_-(y,\gve) \subset \gS_+(y,0)$, we get
\begin{equation}
    \label{Mt16}
    \integdd{\abs{v_\tau\devp{w}{\nu}}}{\gS\backslash\gS_-(y,\gve)}{\gs(x)}{t} \leq \norme{ c_- }_{L^2(\gS)} (m\gve)^{-\frac12}\lrpa{\integdd{e^{-2{\gt_\tau} }\abs{\devp{w}{\nu}}^2 \abs{\devp{\psi}{\nu}}}{\gS_+(y,0)}{\gs(x)}{t}}^{\frac12}.
\end{equation}
Applying the Carleman estimate \eqref{Mt37} and using the fact that $\devp{w}{\nu} = 0$ on $\gS_-(y,0)$, we get that there exist $\tau_1 > 0$ and $C' > 0$ which depend only on $\gO$, $T$, $\norme{q_2}_{L^\infty(Q)}$ and $\norme{A_2}_{L^\infty(Q)^n}$ such that for every $\tau > \tau_1$, we have
\begin{equation}
    \label{Mt17}
    \integdd{e^{-2{\gt_\tau} }\abs{\devp{w}{\nu}}^2 \abs{\devp{\psi}{\nu}}}{\gS_+(y,0)}{\gs(x)}{t} \leq C'\tau^{-1} \norme{e^{-\gt_\tau} Lw }_{L^2(Q)}^2,
\end{equation}
with $L := \devp{}{t} - \gD + A_2\cdot \nabla + q_2$. 
Combining \eqref{Mt15}, \eqref{Mt16} and \eqref{Mt17}, we get
\begin{equation}
    \label{Mt18}
    \forall \tau > \tau_1,\quad \abs{\integdd{v_\tau\devp{w}{\nu}}{\gS\backslash\gS_-(x_0,2\gve)}{\gs(x)}{t}} \leq C'' \tau^{-\frac12} \norme{e^{-\gt_\tau} L w }_{L^2(Q)},
\end{equation}
where $C'' > 0$ depends only on $\gO$, $T$, $A_2$, $q_2$, $\gve$ and $c_-$.
Now by \eqref{Mt3}, we recall that
\begin{equation}
    \label{Mt19}
    e^{-\gt_\tau} L w = (A_2-A_1)\cdot \nabla (c_+ + R_{+,\tau}) + \lrbr{\tau(A_2-A_1) \cdot \nabla \psi + (q_2-q_1)}(c_+ + R_{+,\tau}).
\end{equation}
\new

\subsection{First lemmas:} \label{Flem}

Combining \eqref{Mt11} and \eqref{Mt13}, we obtain 
\begin{equation}
    \label{Mt20}
    \tau I_\tau + II_\tau = - \integdd{v_\tau\devp{w}{\nu}}{\gS\backslash\gS_-(x_0,2\gve)}{\gs(x)}{t}.
\end{equation}
We derive from \eqref{Mt18}, \eqref{Mt19} and \eqref{Mt14} two important intermediate results.
The first one is the previously stated Theorem~\ref{Theo2}. 
The second one, which will be proved after, will be the next step to prove Theorem~\ref{Theo1}.
\new

\begin{Preuve}[of Theorem~\ref{Theo2}]
    We stick to the notations of Section~\ref{Nota} with $A_1 = A_2 =: A$.
    Since $A_1 = A_2$, we obtain $Lw = (q_2-q_1) u_\tau$ from \eqref{Mt12} along with
    \begin{equation*}
        \norme{e^{-\gt_\tau} L w }_{L^2(Q)} = \norme{(q_2-q_1)(c_+ + R_{+,\tau})}_{L^2(Q)} \xrightarrow[\tau\to\infty]{} \norme{(q_2-q_1)c_+}_{L^2(Q)},
    \end{equation*}
    from the decay property \eqref{Mt10}. 
    It follows with \eqref{Mt18} that,
    \begin{equation}
        \label{Mt21}
        \abs{\integdd{v_\tau\devp{w}{\nu}}{\gS\backslash\gS_-(x_0,2\gve)}{\gs(x)}{t}} \xrightarrow[\tau \to \infty]{} 0.
    \end{equation}
    According to \eqref{Mt14} and \eqref{Mt21}, we obtain from \eqref{Mt20},
        \begin{equation}
        \label{Mt22}
        \integdd{(q_2-q_1) c_+c_- }{Q}{x}{t} = 0.
    \end{equation}
    \new

    From now on, without loss of generality, we denote by $q_1$ and $q_2$, their extension by zero over $\inteoo{0}{T} \times \RR^n$.
    We will prove that $q_1 = q_2$ in $ \inteoo{0}{T}\times\RR^n$ using for $t\in\inteoo0T$, the normal operator of the X-ray transform of $(q_2 - q_1)(t,\cdot)$, in the sense of compactly supported distributions over $\RR^n$.
    For this purpose, we will define $h \in L^1_{loc}(\RR^n)$ by
    $$ h : \RR^n\backslash\ense{0} \ni x \longmapsto \frac{1}{\abs{x}^{n-1}} \in \RR^*_+, $$
    and show that for almost every $t \in \inteoo0T$, $(q_2-q_1)(t,\cdot ) \ast h = 0$ on $B(x_0,\gve)$. 
    Applying \cite[Theorem 1.1.]{ilmavirtaUniqueContinuationNormal2020}, we will conclude that $(q_2-q_1)(t,\cdot) = 0$ for almost every $t \in \inteoo{0}{T}$, thus concluding the proof.
    \new

    Let $\chi \in \reguc{\infty}{0,T}$. 
    We consider $\chi_+ = \chi$ and $\chi_-$ such that $\chi_- \equiv 1 $ over $\supp(\chi)$.
    Setting $h_\pm \equiv 1$, we get from \eqref{Mt22} and the definition of the principal parts,
    \begin{equation}
        \label{Mt23}
        \Integd{\lrpa{\integd{(q_2-q_1)(t,x) \abs{x-y}^{-(n-1)} }{\RR^n}{x}}\chi(t)}{0}{T}{t} = 0,
    \end{equation}
    for every $y \in B(x_0,\gve)$.
    Let us consider $\xi \in \reguc{\infty}{B(x_0,\gve)}$. 
    Multiplying \eqref{Mt23} by $\xi(y)$ for every $y \in B(x_0,\gve)$ and integrating over $y \in B(x_0,\gve)$ gives us,
    \begin{equation*}
        \label{Mt24}
        \Integd{\lrbr{\integd{\bbr{ (q_2-q_1)(t,\cdot) \ast h }(y)\xi(y)}{B(x_0,\gve)}{y}}\chi(t)}{0}{T}{t} = 0.
    \end{equation*}
    Recalling that $\chi \in \reguc{\infty}{0,T}$ is arbitrary chosen, we deduce that for almost every $t \in \inteoo0T$,
    \begin{equation*}
        \label{Mt25}
        \integd{\bbr{ (q_2-q_1)(t,\cdot) \ast h }(y)\xi(y)}{B(x_0,\gve)}{y} = 0.
    \end{equation*}
    For such $t \in \inteoo0T$, since $\xi \in \reguc{\infty}{B(x_0,\gve)}$ is arbitrary chosen, we have for almost every $y \in B(x_0,\gve)$,
    \begin{equation*}
        \label{Mt26}
        \lrbr{(q_2-q_1)(t,\cdot) \ast h}(y) = 0.
    \end{equation*}
    In conclusion, we get that for almost every $t \in \inteoo0T$, $(q_2-q_1)(t,\cdot) \ast h = 0$ in $B(x_0,\gve)$ in the sense of compactly supported distributions. 
    Using \cite[Theorem 1.1.]{ilmavirtaUniqueContinuationNormal2020}, we get that $(q_2-q_1)(t,\cdot) = 0$ for almost every $t \in \inteoo0T$ which implies that $q_1 =q_2$. 
\end{Preuve}

The next step of the proof of Theorem~\ref{Theo1} is to prove that the X-ray transform of the zero extension of the continuous vector field $A_2-A_1$ vanishes over every oriented line that intersect $B(x_0,\gve)$.
We will then use \cite[Theorem 1.2.]{ilmavirtaXrayTomographyOneforms2021} to begin the construction of $\gvp$. 
We will come back to the definition of the X-ray transform during the proof of the following lemma,
\begin{Lemme}
    \label{Lem:Mt2}
    Under the hypothesis of Theorem~\ref{Theo1}, for every $t_0 \in \inteoo{0}{T}$, $y \in B(x_0,\gve)$ and $\gt \in \SS^{n-1}$, we get
    \begin{equation}
        \label{Mt27}
        \integd{\lrpa{\frac{A_2 - A_1}{2}}(t_0,y+r\gt)\cdot\gt}{\RR}{r} = 0, 
    \end{equation}
    where $A_2-A_1$ has been continuously extended to $0$ in $\inteoo{0}{T} \times \RR^n \backslash \ade{\gO}$.
\end{Lemme}

\begin{Preuve}
    We stick to the notations of Section \ref{Nota}.
    Using \eqref{Mt18} with \eqref{Mt19}, we get
    \begin{equation}
        \label{Mt28}
        \tau^{-1} \integdd{v_\tau\devp{w}{\nu}}{\gS\backslash\gS_-(x_0,2\gve)}{\gs(x)}{t} \xrightarrow[\tau \to \infty]{} 0.
    \end{equation}
    Dividing \eqref{Mt20} by $\tau$ and using \eqref{Mt14} along with \eqref{Mt28}, we conclude that
    \begin{equation}
        \label{Mt29}
        \integdd{\nabla \psi \cdot (A_2-A_1)c_+c_-}{Q}{x}{t} = 0,
    \end{equation}
    for any $y \in B(x_0,\gve)$ and $c_\pm$ of the form \eqref{Mt6} with any $\tilde{A}_\pm \in \regub{1}{0,T;\regu{0}{\RR^n}} \cap \regub{0}{0,T;\regu{2}{\RR^n}}$ extension of $A_\pm$, $h_\pm \in \reguc{\infty}{\SS^{n-1}}$ and $\chi_\pm \in \reguc{\infty}{0,T}$.
    \new

    Since $\tilde{A}_1$ and $\tilde{A}_2$ are arbitrary chosen, a first approach for this result would be to define $\tilde{A}_2$ from $\tilde{A}_1$, such that $\tilde{A}_2 - \tilde{A}_1$ is the zero extension of $A_2 - A_1$ over $\inteoo{0}{T}\times\RR^n$, that we denote $A$. 
    However, the required regularity for $A_1$ and $A_2$ necessarily implies that ${\nabla A_1}_{|_\gS} = {\nabla A_2}_{|_\gS}$ if we suppose $A = \tilde{A}_2 - \tilde{A}_1$. 
    In order to avoid this condition, we construct the following approximation of $A$; we consider for $m\in\NN$ and $i=1,2$, the extension $\tilde{A}_{i,m} \in \regub{1}{0,T;\regu{0}{\RR^n}} \cap \regub{0}{0,T;\regu{2}{\RR^n}}$ of $A_i$ defined for every $(t,x) \in \inteoo{0}{T}\times \RR^n$ by 
    $$ \tilde{A}_{i,m}(t,x) = \tilde{A}_i(t,x) \rho_m(x), $$
    where $\rho_m \in \regu{\infty}{\RR^n}$ is a smooth cutoff function such that $\rho_m(x) = 1$ for $x \in \ade{\gO}$ and $\rho_m$ is compactly supported in $\gO+B\lrpa{0,\gve m^{-1}}$. 
    By definition, for every $m\in\NN$ and $i=1,2$, we get 
    \begin{equation}
        \label{Mt30}
        \begin{dcases}
            \tilde{A}_{i,m} = A_i &\text{ in } Q \cup \gS, \\
            \tilde{A}_{i,m} = 0 &\text{ in } \inteoo{0}{T} \times \Bbbr{\RR^n \backslash \lrpa{\gO+B\lrpa{0,\gve m^{-1}}}}.
        \end{dcases}
    \end{equation}
    Condition \eqref{Mt30} implies $\tilde{A}_{1,m} = \tilde{A}_{2,m} = 0$ in $\inteoo{0}{T} \times B(x_0,\gve)$ according to \eqref{Hyp1}.
    We notice that $A$ is continuous on $\inteoo{0}{T}\times\RR^n$ since ${A_1}_{|_\gS} = {A_2}_{|_\gS}$ and we also get for every $t_0 \in \inteoo{0}{T}$ that
    \begin{equation}
        \label{Mt31}
        (A_{2,m} - A_{1,m})(t_0,\cdot) \xrightarrow[m\to\infty]{} A(t_0,\cdot),
    \end{equation}
    uniformly over $\RR^n$. 
    \new

    Let us fix $t_0 \in \inteoo{0}{T}$.
    For $F$ a compactly supported, continuous, vector field in $\RR^n$, we define the X-ray transform of $F$ by 
    $$ X_1F(\gg) = \integ{F}{\gg}, $$
    where $\gg$ is an oriented line in $\RR^n$ and $F$ is identified as a $1$-form over $\RR^n$.
    Such oriented line $\gg$ can be described by a couple $(z,\gt) $ in $ \gg \times \SS^{n-1}$ such that
    $$ \RR \ni r \longmapsto z + r\gt \in \gg $$
    is a parametrization of $\gg$. 
    In this case, integration of $1$-form gives us 
    $$ X_1F(\gg) = \integd{F(z+r\gt)\cdot \gt}{\RR}{r}. $$
    The final step of this proof is to compute the X-ray transform of $A(t_0,\cdot)$ over every oriented line $\gg$ in $\RR^n$ that intersect $B(x_0,\gve)$, starting from the integral identity \eqref{Mt29}.
    \new
    
    We introduce the polar coordinates centered at $y$ given by
    $$ \Phi : \inteoo{0}{T} \times \RR_+^* \times \SS^{n-1} \ni (t,r,\gt) \longmapsto (t,y+r\gt) \in \inteoo{0}{T} \times \RR\backslash\ense{y}. $$
    We also consider $\chi \in L_+^1(\RR)$ such that $\supp(\chi) \subset \inteoo{-1}{1}$ and 
    $$\integd{\chi(t)}{\RR}{t} = 1 .$$
    For every $\gd > 0$, we consider $\chi_\gd = \gd \chi(\gd (\cdot-t_0))$ which is supported in $\inteoo{t_0 - \gd^{-1}}{t_0 + \gd^{-1}}$.
    We fix $\gd_0 > 0$ such that for every $\gd > \gd_0$, $\supp(\chi_\gd) \subset \inteoo{0}{T}$.
    For $\gd > \gd_0$, we can use \eqref{Mt29} with the extension $\tilde{A}_{i,m}$ of $A_i$, $\chi_+ = \chi_\gd$, $\chi_-$ such that $\chi_- \equiv 1$ over $\supp(\chi_\gd)$, $h_+ = h$ for $h \in \reguc{\infty}{\SS^{n-1}}$ and $h_- =1$, to get for every $m \in \NN$, that
    \begin{equation}
        \label{Mt32}
        \forall \gd > \gd_0,\quad \Integd{\lrpa{\integd{\lrpa{\ga_m \circ \Phi^{-1}}(t,x) \abs{x-y}^{-(n-1)} h\lrpa{ \frac{x-y}{\abs{x-y}} }}{\gO}{x}}\chi_\gd(t)}{0}{T}{t} = 0,
    \end{equation}
    where $\ga_m : \inteoo{0}{T} \times \RR^*_+ \times \SS^{n-1} \longrightarrow \RR$ is a continuous map defined by
    $$\ga_m(t,r,\gt) = \gt \cdot (\tilde{A}_{2,m} - \tilde{A}_{1,m})(t,y+r\gt) \exp\lrpa{\Integd{\frac{(\tilde{A}_{2,m} - \tilde{A}_{1,m})(t,y+s\gt)}{2}\cdot \gt}{0}{r}{s}}.$$
    Since the map $\inteoo{0}{T} \ni t \longmapsto \integd{\lrpa{\ga_m \circ \Phi^{-1}}(t,x) \abs{x-y}^{-(n-1)} h\lrpa{ \frac{x-y}{\abs{x-y}} }}{\gO}{x}$ is continuous, we conclude by taking the limit $\gd \to \infty$ in \eqref{Mt32} that
    \begin{equation*}
        \label{Mt33}
        \integd{\lrpa{\ga_m \circ \Phi^{-1}}(t_0,x) \abs{x-y}^{-(n-1)} h\lrpa{ \frac{x-y}{\abs{x-y}} }}{\gO}{x} = 0.
    \end{equation*}
    We set $\gO_m := \lrpa{\gO+B\lrpa{0,\gve m^{-1}}} \backslash \ade{\gO}$ for every $m \in\NN$, and we get from the definition of $\tilde{A}_{i,m}$ for $i=1,2$ that
    \begin{equation}
        \label{Mt38}
        \forall m\in \NN,\quad \integd{\lrpa{\ga_m \circ \Phi^{-1}}(t_0,x) \abs{x-y}^{-(n-1)} h\lrpa{ \frac{x-y}{\abs{x-y}} }}{\RR^n \backslash \gO_m}{x} = 0.
    \end{equation}
    Since $\ga_m(t_0,\cdot)$ converge pointwise to the continuous and compactly supported function $\ga(t_0,\cdot)$ defined by
    $$ \forall (r,\gt) \in \RR^*_+ \times \SS^{n-1},\quad \ga(t_0,r,\gt) = \gt \cdot A(t_0,y+r\gt) \exp\lrpa{\Integd{\frac{A(t_0,y+s\gt)}{2}\cdot \gt}{0}{r}{s}}, $$
    and $\indi{\gO_m} \xrightarrow[m \to \infty]{} 0$ pointwise, then by applying Lebesgue dominated convergence theorem and taking the limit $m\to \infty$ of \eqref{Mt38}, we get 
    \begin{equation*}
        \integd{\lrpa{\ga \circ \Phi^{-1}}(t_0,x) \abs{x-y}^{-(n-1)} h\lrpa{ \frac{x-y}{\abs{x-y}} }}{\RR^n}{x} = 0.
    \end{equation*}
    After a change of coordinates, the latest equality gives us in polar coordinates centered at $y$, 
    \begin{equation}
        \label{Mt34}
        \integd{\lrpa{\Integd{ \ga(t_0,r,\gt) }{0}{+\infty}{r} } h(\gt)}{\SS^{n-1}}{\gt} = 0.
    \end{equation}
    Recalling that $h \in \reguc{\infty}{\SS^{n-1}}$ is arbitrary chosen and the fact that $\SS^{n-1} \ni \gt \longmapsto \Integd{ \ga(t_0,r,\gt) }{0}{+\infty}{r} $ is continuous, we deduce from \eqref{Mt34} that for every $\gt \in \SS^{n-1}$,
    \begin{equation*}
        \label{Mt35}
        \Integd{ \ga(t_0,r,\gt) }{0}{+\infty}{r} = 0.
    \end{equation*}
    Using fundamental theorem of calculus, we get that
    $$ \exp\lrpa{\Integd{\frac{A(t_0,y+s\gt)}{2}\cdot \gt}{0}{+\infty}{s}} - 1 = 0, $$
    which gives us for every $\gt \in \SS^{n-1}$,
    \begin{equation}
        \label{Mt36}
        \Integd{\frac{A(t_0,y+s\gt)}{2}\cdot \gt}{0}{+\infty}{s} = 0.
    \end{equation}
    In order to conclude, we consider \eqref{Mt36} with $-\gt$ instead of $\gt$ and the change of variables $u=-s$ to get \eqref{Mt27}.
\end{Preuve}

\subsection{Completion of the proof of Theorem~\ref{Theo1}:}

Now that we know that for a fixed $t_0$ in $\inteoo{0}{T}$, the X-ray transform of $\frac{A_2-A_1}{2}(t_0,\cdot)$ extended by zero over $\RR^n\setminus\ade{\gO}$ vanishes over every oriented line that intersect $B(x_0,\gve)$, we can use \cite[Theorem 1.2.]{ilmavirtaXrayTomographyOneforms2021} with Proposition~\ref{Pro:Ca4} in order to prove the following,
\begin{Proposition}
    \label{Pro:Mt1}
    Under the hypothesis of Theorem~\ref{Theo1}, there exists $\gvp \in \gvpSpace$ such that
    $$ A_2-A_1 = 2\nabla \gvp \quad \text{ in } Q .$$
\end{Proposition}

\begin{Preuve}
    Let us note $A$, the continuous zero extension of $A_2-A_1$ over $\inteoo{0}{T} \times \RR^n\setminus\ade{\gO}$ and let us define $\gvp : \inteoo{0}{T} \times \RR^n \rightarrow \RR$ with 
    $$ \forall (t,x) \in \inteoo{0}{T} \times \RR^n, \quad \gvp(t,x) = \frac{1}{2}\Integd{A(t,x_0 + s(x-x_0))\cdot (x-x_0)}{0}{1}{s}. $$
    We set $t_0 \in \inteoo{0}{T}$. 
    From Lemma~\ref{Lem:Mt2}, we know that the X-ray transform of $A(t_0,\cdot)$ vanishes over every oriented line that intersect $B(x_0,\gve)$.
    According to \cite[Theorem 1.2.]{ilmavirtaXrayTomographyOneforms2021}, we get $\dif{A(t_0,\cdot)}{}{} =0$ over $\RR^n$ where $\td$ denotes the exterior derivative operator for $1$-forms in the sense of distributions.
    Applying Proposition~\ref{Pro:Ca4} (see the Appendix \ref{ConsGaug}), we obtain that, $\gvp(t_0,\cdot)$ belongs to $\regu{1}{\RR^n}$ with $\nabla \gvp(t_0,\cdot) = \frac{A}{2}(t_0,\cdot)$ which implies
    \begin{equation*}
        A_2-A_1 = 2\nabla \gvp \quad \text{ in } Q .
    \end{equation*}
    Since $\RR^n \backslash \gO$ is connected and $\gvp(\cdot, x_0) = 0$ with $\nabla \gvp = 0$ over $\inteoo{0}{T} \times \RR^n \backslash \ade{\gO}$, we get that $\gvp_{|_{\gS}} = 0$ by continuity of $\gvp$. 
    Using also the fact $A_2-A_1 = 2\nabla \gvp$ over $\inteoo{0}{T} \times \RR^n$, we get that $\devp{\gvp}{\nu}_{|_{\gS}} = 0$ since ${A_1}_{|_\gS} = {A_2}_{|_\gS}$.
    \new

    We now prove that $\gvp \in \gvpSpace$. 
    For any compact set $K'$ of $\RR^n$, let $\norme{\cdot}_{K',\infty}$ be defined for every $F \in L_{loc}^\infty(\RR^n)^n$ by $\norme{F}_{K',\infty}= \norme{F\indi{K'}}_{L^\infty(\RR^n)^n}$. 
    Recall that $A \in  \regub{0}{0,T; \regu{0}{\RR^n}}^n$ and for every $t\in\inteoo0T$, $A(t,\cdot)$ is supported over $\ade{\gO}$. 
    Moreover, for every $t$, $t' \in \inteoo{0}{T}$ and compact set $K$ of $\RR^n$, we obtain
    \begin{equation*}
        \norme{A(t,\cdot) - A(t',\cdot)}_{K,\infty} \leq \norme{A_2(t,\cdot) - A_2(t',\cdot)}_{K\cap \ade{\gO},\infty} + \norme{A_1(t,\cdot) - A_1(t',\cdot)}_{K\cap \ade{\gO},\infty}.
    \end{equation*}
    In view of Proposition \ref{Pro:Ca5}, we have $\gvp \in \regub{0}{0,T; \regu{0}{\RR^n}}$ and $\nabla \gvp = A \in \regub{0}{0,T; \regu{0}{\RR^n}}^n$ which implies that $\gvp \in \regub{0}{0,T; \regu{1}{\ade{\gO}}}$. 
    In addition, for every $t \in \inteoo{0}{T}$, we find
    $$ \nabla \gvp(t,\cdot ) = \frac{A_2(t,\cdot)- A_1(t,\cdot)}{2} \quad \text{ in } \gO, $$
    with $A_2$ and $A_1$ in $\regub{0}{0,T; \regub{2}{\gO}}^n$ and it follows that $\gvp \in \regub{0}{0,T; \regub{2}{\gO}}$. 
    Thus, $\gvp$ is lying in $\regub{0}{0,T; \regub{2}{\gO} \cap \regu{1}{\ade{\gO}}}$. 
    Since $A \in \regub{1}{0,T; \regu{0}{\RR^n}}^n$, we get that 
    $$ \forall (t,x) \in \inteoo{0}{T} \times \RR^n,\quad \devp{\gvp}{t}(t,x) = \frac12\Integd{\devp{A}{t}(t,x_0 + s(x-x_0))\cdot (x-x_0)}{0}{1}{s}. $$
    Applying Proposition~\ref{Pro:Ca5} with the fact that $\devp{A}{t}$ belongs to $\regub{0}{0,T; \regu{0}{\RR^n}}^n$, we conclude that $\gvp \in \regub{1}{0,T; \regu{0}{\RR^n}}$ and thus $\gvp \in \gvpSpace$.
\end{Preuve}

We are now able to conclude the proof of Theorem~\ref{Theo1}.
To summarize, for $(A_1, q_1)$ and $(A_2, q_2)$ in $\ASpace \times L^\infty(Q)$, if $\gL_{(A_1, q_1)} = \gL_{(A_2, q_2)}$ and ${A_1}_{|_\gS} = {A_2}_{|_\gS}$ then there exists $\gvp \in G$ such that $A_2 - A_1 = 2 \nabla \gvp$ thanks to Proposition~\ref{Pro:Mt1}. 
Using the gauge invariance \eqref{Eq:InvGau} we also know that 
$$ \gL_{(A_2, q_2)} = \gL_{S_{-\gvp}(A_1 + 2\nabla \gvp, q_2)} = \gL_{(A_1, q_2 + \devp{\gvp}{t} - \gD \gvp + \abs{\nabla \gvp}^2 + A_1\cdot \nabla \gvp )} $$
which gives $\gL_{(A_1, q_1)} = \gL_{(A_1, q_2 + \devp{\gvp}{t} - \gD \gvp + \abs{\nabla \gvp}^2 + A_1\cdot \nabla \gvp )}$.
Applying Theorem~\ref{Theo2}, we get that
$$ q_1 = q_2 + \lrpa{\devp{\gvp}{t} - \gD \gvp + \abs{\nabla \gvp}^2 + A_1\cdot \nabla \gvp}, $$
thus giving us $(A_2,q_2) = S_\gvp(A_1,q_1)$ and concluding the proof of Theorem~\ref{Theo1}.

\section{Global uniqueness:}   \label{FuCo}

This section is devoted to problem {\bf (IP2)}. 
We already stated in Theorem~\ref{Theo2}, that if we fix $A \in \ASpace$, then for every $q_1$, $q_2 \in \qSpace$, we have
    $$ \gL_{(A,q_1)} = \gL_{(A,q_2)} \implies q_1 = q_2 .$$
We would like to extend this property for particular cases of pair $(A,q)$ in order to obtain the following implication 
$$ \gL_{(A_1,q_1)} = \gL_{(A_2,q_2)} \implies (A_1,q_1) = (A_2,q_2) .$$

We recall the following classic results,
\begin{Lemme}
    \label{Lem1}
    Let $V \in L^\infty(Q)^n$. If $u \in H^{1,2}(Q)$ satisfies the following conditions,
    \begin{equation}
        \label{F1}
        \begin{dcases}
            \devp{u}{t} - \gD u + V \cdot \nabla u = 0 , &\quad \text{ in } Q,\\
            u_{|_\gS} = 0, &\quad \text{ on } \gS, \\
            \devp{u}{\nu}_{|_\gS} = 0, &\quad \text{ on } \gS,
        \end{dcases}
    \end{equation}
    then $u$ vanishes in $Q$.
\end{Lemme}

\begin{Preuve}
    The proof of this lemma relies mainly on the unique continuation property of the linear parabolic equation stated above. 
    We give the proof of this lemma for the sake of completeness.
    Let us consider $u$ a solution of \eqref{F1}, a connected open set $\gO'$ of $\RR^n$ such that $\ade{\gO} \subset \gO'$, and $\tilde{V}$ (resp. $\tilde{u}$) the zero extension over $Q' := \inteoo{0}{T}\times \gO'$ of $V$ (resp. $u$). 
    We have $\tilde{V} \in L^\infty(Q')^n$ and, since $u \in H^{1}\lrpa{ 0,T;L^2(\gO) }$, we get that $\tilde{u} \in H^{1}\lrpa{0,T;L^2(\gO')}$. 
    In addition, since $\devp{u}{\nu}_{|_\gS} = 0$ and $u_{|_\gS} = 0$, we obtain that $\nabla u_{|_\gS} = 0$. 
    Using the fact that $\gO$ has a $\cC^1$ boundary, it follows that $\tilde{u}(t,\cdot) \in H^1_0(\gO')$ for almost every $t \in \inteoo{0}{T}$. 
    We then use the $\cC^2$ boundary of $\gO$ along with the elliptic regularity, and the fact that $u \in L^2(0,T;H^2(\gO) \cap H^1_0(\gO))$ to conclude that $\tilde{u} \in L^2(0,T;H^2(\gO') \cap H^1_0(\gO'))$ and verify 
    \begin{equation*}
        \label{F2}
        \devp{\tilde{u}}{t} - \gD \tilde{u} + \tilde{V} \cdot \nabla \tilde{u} = 0 , \quad \text{ in } Q'.
    \end{equation*}
    By setting $(y,r) \in \gO' \times \RR^*_+$ such that $\ade{B(y,r)} \cap \ade{\gO} = \emptyset$ and $B(y,r) \subset \gO'$, we also find that $\tilde{u}$ vanishes over $\inteoo{0}{T} \times B(y,r)$. 
    Applying unique continuation property of the parabolic operator $\devp{}{t} - \gD + \tilde{V} \cdot \nabla $ \cite[Corollary 1.2.]{sautUniqueContinuationEvolution1987}, we get that $\tilde{u}$ vanishes over $Q'$ thus concluding the proof.
\end{Preuve}

\subsection{When the zeroth order coefficient is known:}

For $q$ known, we consider the problem of determination of $A$ from $\gL_{(A,q)}$.
To do so, we can use unique continuation property of parabolic equation and properties of the gauge function to prove the following, 

\begin{Preuve}[of Corollary~\ref{Coro1}]
    Let us fix $q \in L^\infty(Q)$ along with $A_1$, $A_2$ in $\ASpace$ such that ${A_1}_{|_\gS} = {A_2}_{|_\gS}$ and $\gL_{(A_1,q)} = \gL_{(A_2,q)}$. 
    Following Theorem~\ref{Theo1}, there exists $\gvp \in \gvpSpace$ such that 
    $$ (A_2,q) = \bpa{A_1 + 2\nabla \gvp, q - \lrpa{ \devp{\gvp}{t} - \gD \gvp + (\nabla \gvp + A_1)\cdot \nabla \gvp }} .$$
    We deduce that $\gvp \in \SolSpace$ satisfies 
    \begin{equation*}
        \label{F3}
        \begin{dcases}
            \devp{\gvp}{t} - \gD \gvp + (\nabla \gvp + A_1) \cdot \nabla \gvp = 0 , &\quad \text{ in } Q,\\
            \gvp_{|_\gS} = 0, &\quad \text{ on } \gS, \\
            \devp{\gvp}{\nu}_{|_\gS} = 0, &\quad \text{ on } \gS.
        \end{dcases}
    \end{equation*}
    Considering $V = \nabla \gvp + A_1 \in L^\infty(Q)^n$, we use Lemma~\ref{Lem1} to conclude that $\gvp = 0$ in $Q$ thus giving $A_2 = A_1$.
\end{Preuve}

\subsection{The Fokker-Planck equation:}

Another answer of {\bf (IP2)} is the uniqueness for the partial Dirichlet-to-Neumann operator associated with the Fokker-Planck equation. 
Let us consider for $A \in L^\infty(Q)^n$ such that $\dive(A) \in L^\infty(Q)$, the following evolution equation, namely the Fokker-Planck equation,
\begin{equation}
    \label{PFK2}
    \begin{dcases}
        \devp{u}{t} - \gD u + \dive\lrpa{Au} = 0 , &\quad \text{ in } Q,\\
        u_{|_\gS} = g, &\quad \text{ on } \gS, \\
        u_{|_{t=0}} = 0, &\quad \text{ on } \gO.
    \end{dcases}
\end{equation}
We associate with \eqref{PFK2}, the following partial DN map
$$ \cF_A : \hdirz \ni g \longmapsto \devp{u}{\nu}_{|_{\gS_-(x_0,2\gve)}} \in H^{\frac14,\frac12}(\gS_-(x_0,2\gve)), $$
where $u$ is the unique solution of \eqref{PFK2}.
One can check that
\begin{equation}
    \label{F4}
    \forall A \in \ASpace,\quad \cF_A = \gL_{(A,\dive(A))} .
\end{equation}
\new

\begin{Preuve}[of Corollary~\ref{Coro2}]
    Let us consider $A_1$, $A_2$ in $\ASpace$ such that ${A_1}_{|_\gS} = {A_2}_{|_\gS}$ and assume that we have $\cF_{A_1} = \cF_{A_2}$. 
    Using \eqref{F4} we get that $\gL_{(A_1, \dive(A_1))} = \gL_{(A_2, \dive(A_2))}$. 
    Applying Theorem~\ref{Theo1}, we know that there exists $\gvp \in \gvpSpace$ such that 
    $$ (A_2,\dive(A_2)) = (A_1 + 2\nabla \gvp, \dive\lrpa{A_1} - \lrpa{ \devp{\gvp}{t} - \gD \gvp + (\nabla \gvp + A_1)\cdot \nabla \gvp }) .$$
    Since $A_2 = A_1 + 2 \nabla \gvp$, we get that $\dive\lrpa{ A_2 } = \dive\lrpa{A_1} + 2\gD \gvp$ in $Q$.
    It follows that
    $$ -\devp{\gvp}{t} -\gD \gvp -(\nabla \gvp + A_1)\cdot \nabla \gvp = 0 \text{ in } Q. $$
    Thus, fixing $\ade{\gvp}$ and $\ade{A_1}$ defined by $\ade{\gvp}(t,x) = \gvp(T-t,x)$ and $\ade{A_1}(t,x) = A_1(T-t,x)$, for $(t,x) \in Q$, we get
    \begin{equation*}
        \label{F5}
        \begin{dcases}
            \devp{\ade{\gvp}}{t} - \gD \ade{\gvp} + (\nabla \ade{\gvp} + \ade{A_1}) \cdot \nabla \ade{\gvp} = 0 , &\quad \text{ in } Q,\\
            \ade{\gvp}_{|_\gS} = 0, &\quad \text{ on } \gS, \\
            \devp{\ade{\gvp}}{\nu}_{|_\gS} = 0, &\quad \text{ on } \gS.
        \end{dcases}
    \end{equation*}
    Using Lemma~\ref{Lem1}, we deduce that $\ade{\gvp}$ vanishes in $Q$ which implies $A_2 = A_1$.
\end{Preuve}

\appendix
\section{Auxiliary results:}   \label{Auxi}

This section is devoted to the proof of auxiliary results used throughout the proof of Theorem~\ref{Theo1}.
In all this section, $\anorme_2$ (resp. $\anorme_\infty$) will denotes the usual norm of $L^2(Q)$, $L^2(Q)^n$ and $L^2(\gO)$ depending on the situation (resp. $L^\infty(Q)$ and $L^\infty(Q)^n$).

\subsection{Carleman estimate:}            \label{Carleman}

\newcommand{\vL}{u}
\newcommand{\vPt}{v}
\newcommand{\vPtl}{w}
\newcommand{\para}{\gl}
For this subsection, we will consider $A$ in $L^\infty(Q)^n$ and $q$ in $L^\infty(Q)$ along with $\psi \in \regu{1}{\ade{\gO}} \cap \regub{2}{\gO}$  satisfying the eikonal equation \eqref{Hyp2} and \eqref{Hyp3}. 
For such $\psi$, we define the sets 
$$ \begin{dcases}
    \gG_{\pm,\psi} := \ense{ x \in \gG \given \pm \devp{\psi}{\nu}(x) > 0} ,\\
    \gS_{\pm,\psi} := \inteoo{0}{T}\times \gG_{\pm,\psi}.
\end{dcases} $$
We also recall that the definition of our weight function $\gt_\tau$ is
$$ \forall (t,x) \in \ade{Q},\quad \gt_\tau(t,x) = \tau^2 t + \tau \psi(x) .$$
We fix $ L_\pm = \pm \devp{}{t} - \gD \pm A \cdot \nabla + q $ and derive the following Carleman estimate:
\begin{Theoreme}
    \label{Theo:Ca1}
    For every $\vL_\pm \in H^{1,2}(Q)$ satisfying ${\vL_\pm}_{|_\gS} = 0$, ${\vL_+}_{|_{t=0}} = 0$ and ${\vL_-}_{|_{t=T}} = 0$, there exist $\tau_0 > 0$ and $C >0$ which depend only on $\gO$, $T$, $\psi$, $\norme{q}_\infty$ and $\norme{A}_\infty$ such as, for every $\tau > \tau_0$, we have
    \begin{multline}
        \label{ICa1}
        \tau \integdd{\abs{e^{\mp\gt_\tau}\devp{\vL_\pm}{\nu}}^2 \abs{\devp{\psi}{\nu}}}{\gS_{\pm,\psi}}{\gs(x)}{t} + \tau^2 \integdd{\abs{e^{\mp\gt_\tau}\vL_\pm}^2}{Q}{x}{t} + \tau \integdd{\abs{\nabla\lrpa{e^{\mp\gt_\tau}\vL_\pm}}^2}{Q}{x}{t} \\ \leq C\lrbr{ \integdd{\abs{e^{\mp\gt_\tau} L_\pm \vL_\pm}^2}{Q}{x}{t} + \tau \integdd{\abs{e^{\mp\gt_\tau}\devp{\vL_\pm}{\nu}}^2 \abs{\devp{\psi}{\nu}}}{\gS_{\mp,\psi}}{\gs(x)}{t} }.
    \end{multline}
\end{Theoreme}

\begin{remark}
    When $\psi : x \in \ade{\gO} \longmapsto \abs{x-y} \in \RR$ for $y \in \RR^n\setminus\ade{\gO}$, one can check that $\gG_{\pm,\psi} = \gG_\pm(y,0)$ and derive Theorem~\ref{The:ICA}.
\end{remark}

In view of Theorem~\ref{Theo:Ca1}, weights function of the form \eqref{Go16} with $\psi$ satisfying \eqref{Hyp2} and \eqref{Hyp3} can be seen as the parabolic analogue of the concept of a limiting Carleman weight for elliptic equations (see, for example, \cite{bukhgeimRecoveringPotentialPartial2002}, \cite{dossantosferreiraLimitingCarlemanWeights2009} and \cite{krupchykUniquenessInverseBoundary2014}).
As far as we know, so far, such parabolic limiting Carleman weight were limited to linear functions (see \cite{choulliLogarithmicStabilityDetermining2018}, \cite{kianDeterminingBreakingGauge2024}).

\begin{Preuve}[of Theorem~\ref{Theo:Ca1}]
    In all the remaining parts of the proof, we set $T_+ = T$ and $T_- = 0$.
    Let us fix $\vL_\pm \in \SolSpace$ such as ${\vL_\pm}_{|_\gS} = 0$, ${\vL_+}_{|_{t=0}} = 0$, ${\vL_-}_{|_{t=T}} = 0$. 
    We consider $ \vPt_\pm =e^{\mp\gt_\tau} \vL_\pm$ on $\ade{Q}$.
    We recall that for every $\tau >0$, we have
    $$ e^{\mp\gt_\tau}(\pm\devp{}{t} - \gD \pm A\cdot \nabla)\vL_\pm = P_{\pm,\tau} \vPt_\pm , $$
    where $P_{\pm,\tau} = \mp \devp{}{t} - \gD \pm (A - 2\tau \nabla \psi) \cdot \nabla + (\tau A\cdot \nabla \psi \mp \tau \gD \psi)$.
    In order to complete the proof of Theorem~\ref{Theo:Ca1}, it suffices to show that there exist $\tau_0 >0$ and $C>0$ which depend only on $\gO$, $T$, $\psi$ and $\norme{A}_\infty$ such as
    \begin{align}
        \label{ICa2}
        \tau \integdd{\abs{\devp{\vPt_\pm}{\nu}}^2 \abs{\devp{\psi}{\nu}}}{\gS_{\pm,\psi}}{\gs(x)}{t} &+ \tau^2 \norme{\vPt_\pm}^2_2 + \tau \norme{\nabla \vPt_\pm}^2_2 \nonumber \\
        &\leq C\lrbr{ \norme{P_{\pm,\tau} \vPt_\pm}^2_2 + \tau \integdd{\abs{\devp{\vPt_\pm}{\nu}}^2 \abs{\devp{\psi}{\nu}}}{\gS_{\mp,\psi}}{\gs(x)}{t} }.
    \end{align}
    That last statement comes from the fact that 
    $$ \norme{e^{\mp\gt_\tau}L_\pm \vL_\pm}^2_2 = \norme{P_{\pm,\tau}\vPt_\pm + q \vPt_\pm}^2_2 \geq \frac{\norme{P_{\pm,\tau}\vPt_\pm}^2_2}{2} - \norme{q}^2_\infty \norme{\vPt_\pm}^2_2. $$
    We will then demonstrates that \eqref{ICa2} holds true for $\tau$ sufficiently large and $C$ depending only on $\gO$, $T$, $\psi$ and $\norme{A}_\infty$.
    \new

    Known proof of an analogous Carleman estimate in the elliptic case requires to convexify the weight. 
    Unlike the elliptic case, we can prove our result by fixing a parameter $\gl >0$ that will be chosen large enough, and setting
    $$ \tilde{\gt}_{\tau,\gl}(t,x) := \gt_\tau(t,x) - \gl\psi(x) = \tau^2 t + (\tau - \gl)\psi(x),\quad \forall (t,x) \in \ade{Q}, $$
    along with $\vPtl_\pm := e^{\mp\tilde\gt_{\tau,\gl}}\vL_\pm = e^{\pm \gl\psi}\vPt_\pm$ and
    $$ e^{\mp\tilde\gt_{\tau,\gl}}L_\pm \vL_\pm = P^A_{\pm,\tau,\gl} \vPtl_\pm, $$
    with $P^A_{\pm,\tau,\gl} = P^0_{\pm,\tau,\gl} \pm A\cdot\nabla + (\tau-\gl)A\cdot\nabla\psi$, where $P^0_{\pm,\tau,\gl} = P^{(1)}_{\pm,\tau,\gl} + P^{(2)}_{\pm,\tau,\gl} $ and
    $$ \begin{dcases}
    P^{(1)}_{\pm,\tau,\gl} := 2\gl\tau + [-\gD - \gl^2 \pm (\gl-\tau)\gD\psi] ,\\
    P^{(2)}_{\pm,\tau,\gl} := \pm \lrbr{\devp{}{t} + 2(\gl - \tau)\nabla \psi\cdot\nabla}.
        \end{dcases} $$
    \new

    In order to prove \eqref{ICa2}, we will give a suitable lower bound for $\norme{P^A_{\pm,\tau,\gl}\vPtl_\pm}^2_2$.
    We begin with the following inequality:
    \begin{align}
        \label{ICa3}
        \norme{P^A_{\pm,\tau,\gl}\vPtl_\pm}^2_2 &\geq \frac{\norme{P^0_{\pm,\tau,\gl}\vPtl_\pm}^2_2}{2} - \norme{(\tau - \gl)(A\cdot\nabla \psi)\vPtl_\pm \pm A\cdot \nabla \vPtl_\pm}^2_2 \nonumber \\
        &\geq \frac{\norme{P^0_{\pm,\tau,\gl}\vPtl_\pm}^2_2}{2} - C_2 \norme{A}^2_\infty \lrbr{(\tau^2+\gl^2)\norme{\vPtl_\pm}^2_2 + \norme{\nabla \vPtl_\pm}^2_2} .
    \end{align}
    In all this proof, $C_i$ ($i =1,2$) denotes a generic positive constant depending on $\gO$, $T$ and $\psi$ that might change from line to line.
    \new 
    
    Let us now give a lower bound of $\norme{P^0_{\pm,\tau,\gl} \vPtl_\pm}^2_2$.
    We recall that
    \begin{equation}
        \label{ICa4}
        \norme{P^0_{\pm,\tau,\gl} \vPtl_\pm}^2_2 \geq \norme{P^{(1)}_{\pm,\tau,\gl}\vPtl_\pm}^2_2 + 2 \integdd{\lrpa{P^{(1)}_{\pm,\tau,\gl}\vPtl_\pm} \lrpa{P^{(2)}_{\pm,\tau,\gl}\vPtl_\pm}}{Q}{x}{t}.
    \end{equation}
    Since ${\vPtl}_{\pm|_\gS} = 0$, we get
    \begin{align}
        \label{ICa5}
        \norme{P^{(1)}_{\pm,\tau,\gl}\vPtl_\pm}^2_2 =& \integdd{\abs{2 \gl \tau\vPtl_\pm + [-\gD - \gl^2 \pm (\gl-\tau)\gD\psi]\vPtl_\pm}^2}{Q}{x}{t} \nonumber \\
        \geq& \integdd{\abs{2\gl \tau\vPtl_\pm}^2}{Q}{x}{t} + 2 \integdd{\bpa{2\gl \tau\vPtl_\pm}\bpa{-\gD\vPtl_\pm }}{Q}{x}{t} \nonumber \\
        &+ 2 \integdd{\bpa{2\gl \tau\vPtl_\pm}\bpa{- \gl^2 \pm (\gl-\tau)\gD\psi}\vPtl_\pm}{Q}{x}{t} \nonumber \\
        \geq& 4\gl^2\tau^2 \norme{\vPtl_\pm}^2_2 + 4\gl \tau\integdd{\abs{\nabla \vPtl_\pm}^2}{Q}{x}{t} + 4\gl \tau \integdd{[\pm (\gl - \tau)\gD\psi - \gl^2] \vPtl_\pm^2}{Q}{x}{t} \nonumber \\
        \geq& C_1\lrbr{\gl^2\tau^2 \norme{\vPtl_\pm}^2_2 + \gl\tau \norme{\nabla \vPtl_\pm}^2_2} - C_2\lrbr{\gl\tau^2 + \gl^3\tau + \gl^2\tau}\norme{\vPtl_\pm}^2_2.
    \end{align}
    We write
    \begin{equation*}
        2 \integdd{\lrpa{P^{(1)}_{\pm,\tau,\gl}\vPtl_\pm} \lrpa{P^{(2)}_{\pm,\tau,\gl}\vPtl_\pm}}{Q}{x}{t} = I_\pm +II_\pm + III_\pm + IV_\pm,
    \end{equation*}
    where 
    \begin{align*}
        I_\pm &= \pm 2\integdd{\lrbr{2\gl \tau - \gl^2 \pm (\gl - \tau)\gD\psi}\vPtl_\pm \devp{\vPtl_\pm}{t}}{Q}{x}{t}, \\ 
        II_\pm &= \mp 2\integdd{\gD\vPtl_\pm \devp{\vPtl_\pm}{t}}{Q}{x}{t}, \quad
        III_\pm = \pm 4(\tau-\gl)\integdd{\gD\vPtl_\pm \nabla \psi \cdot \nabla \vPtl_\pm}{Q}{x}{t}, \\
        IV_\pm &= \pm 4(\gl-\tau)\integdd{\lrbr{2\gl \tau - \gl^2 \pm (\gl-\tau) \gD\psi} \vPtl_\pm \nabla \psi \cdot \nabla \vPtl_\pm}{Q}{x}{t}.
    \end{align*}
    \new

    Since ${\vPtl}_{\pm|_\gS} = 0$, we have $\nabla\vPtl_{\pm|_\gS} = \lrpa{\devp{\vPtl_\pm}{\nu}} \nu$.
    Noting $K_\pm = 2\gl \tau - \gl^2 \pm (\gl -\tau)\gD \psi$, and recalling that ${\vPtl}_{+|_{t=0}}=0$ and ${\vPtl}_{-|_{t=T}}=0$, we get
    $$ I_\pm = 2\integdd{K_\pm\vPtl_\pm \devp{\vPtl_\pm}{t}}{Q}{x}{t} = \integd{K_\pm\vPtl_\pm^2(T_\pm,x)}{\gO}{x},$$
    which implies
    \begin{equation}
        \label{ICa6}
        I_\pm \geq \lrbr{C_1 \gl\tau - C_2 (\tau + \gl^2+ \gl)}\norme{\vPtl_\pm(T_\pm,\cdot)}^2_2.
    \end{equation}

    For $II_\pm$, we integrate by parts to get 
    \begin{equation}
        \label{ICa7}
        II_\pm = 2 \integdd{\devp{\nabla \vPtl_\pm}{t} \cdot \nabla\vPtl_\pm}{Q}{x}{t} = \integdd{\devp{\lrbr{\abs{\nabla \vPtl_\pm}^2}}{t}}{Q}{x}{t} = \integd{\abs{\nabla \vPtl_\pm(T_\pm,x)}^2}{\gO}{x} \geq 0.
    \end{equation}

    As for $III_\pm$ we can see that, by integrating by parts, we find
    \begin{align*}
        \integdd{\gD\vPtl_\pm \nabla \psi \cdot \nabla \vPtl_\pm}{Q}{x}{t} =& \integdd{(\devp{\vPtl_\pm}{\nu})\nabla\psi\cdot\nabla\vPtl_\pm}{\gS}{\gs(x)}{t} - \integdd{\nabla\vPtl_\pm \cdot \nabla\lrpa{\nabla \vPtl_\pm \cdot \nabla \psi}}{Q}{x}{t} \\
        =& \integdd{(\devp{\vPtl_\pm}{\nu})^2\devp{\psi}{\nu}}{\gS}{\gs(x)}{t} - \integdd{ \Dif{\psi}{}{2}(\nabla\vPtl_\pm, \nabla\vPtl_\pm) }{Q}{x}{t} \\
        & - \frac12 \integdd{ \nabla\psi \cdot \nabla\lrpa{\abs{\nabla\vPtl_\pm}^2} }{Q}{x}{t} \\
        =& \frac12 \integdd{(\devp{\vPtl_\pm}{\nu})^2\devp{\psi}{\nu}}{\gS}{\gs(x)}{t} + \frac12 \integdd{\abs{\nabla\vPtl_\pm}^2 \gD\psi}{Q}{x}{t} \\
        & - \integdd{ \Dif{\psi}{}{2}(\nabla\vPtl_\pm, \nabla\vPtl_\pm) }{Q}{x}{t}.
    \end{align*}
    Recalling that $\psi \in \regub{2}{\gO}$ and that on $\gG$, we have $\devp{\psi}{\nu} = \indi{\gG_+(y,0)}\abs{\devp{\psi}{\nu}} - \indi{\gG_-(y,0)}\abs{\devp{\psi}{\nu}}$, we obtain
    \begin{align}
        \label{ICa8}
        III_\pm \geq & 2(\tau-\gl) \integdd{(\devp{\vPtl_\pm}{\nu})^2 \devp{\psi}{\nu}}{\gS_{\pm,\psi}}{\gs(x)}{t} + 2(\gl-\tau) \integdd{(\devp{\vPtl_\pm}{\nu})^2 \devp{\psi}{\nu}}{\gS_{\mp,\psi}}{\gs(x)}{t} \nonumber \\ 
        & - C_2(\tau+\gl)\norme{\nabla\vPtl_\pm}^2_2.
    \end{align}

    For $IV_\pm$, using \eqref{Hyp3}, we notice for $k=0,1$ that
    \begin{equation*}
        \integdd{(\gD\psi)^k \nabla\psi \cdot \nabla \vPtl_\pm^2}{Q}{x}{t} = -\integdd{\vPtl_\pm^2 \dive\lrpa{(\gD\psi)^k \nabla \psi}}{Q}{x}{t} \geq -C_2\norme{\vPtl_\pm}_2^2,
    \end{equation*}
    from which we obtain that 
    \begin{align}
        IV_\pm &= 2(\gl-\tau)\integdd{\lrbr{ \pm 2\gl \tau \mp \gl^2 + (\gl-\tau) \gD\psi}\nabla \psi \cdot \nabla \vPtl_\pm^2}{Q}{x}{t}  \nonumber \\
        \label{ICa9} &\geq -C_2 \lrbr{\gl\tau^2 + \tau^2 + \gl^2\tau + \gl\tau + \gl^3 + \gl^2 } \norme{\vPtl_\pm}^2_2.
    \end{align}
    \new

    Combining \eqref{ICa5}-\eqref{ICa9} with \eqref{ICa4} and accordingly to \eqref{ICa3}, we obtain
    \begin{align}
        \label{ICa10}
        \norme{P^A_{\pm,\tau,\gl} \vPtl_\pm}^2_2 \geq 
        &\lrbr{C_1\gl^2\tau^2 - C_2\lrpa{\bpa{\gl + 1 + \norme{A}^2_\infty}\tau^2 + \lrpa{\gl^3 + \gl^2 + \gl}\tau + \gl^3 + \gl^2 + \norme{A}^2_\infty \gl^2 }}\norme{\vPtl_\pm}^2_2 \nonumber \\
        &+ \lrbr{C_1 \gl \tau - C_2\lrpa{\tau + \gl + \norme{A}^2_\infty}}\norme{\nabla\vPtl_\pm}^2_2 \nonumber \\
        &+ \lrbr{C_1\tau - C_2 \gl}\integdd{(\devp{\vPtl_\pm}{\nu})^2\abs{\devp{\psi}{\nu}}}{\gS_{\pm,\psi}}{\gs(x)}{t} \\
        &+ \lrbr{C_1 \gl - C_2\tau}\integdd{(\devp{\vPtl_\pm}{\nu})^2\abs{\devp{\psi}{\nu}}}{\gS_{\mp,\psi}}{\gs(x)}{t} \nonumber \\
        &+ \lrbr{C_1\gl \tau - C_2\lrpa{\tau+\gl^2+\gl}}\norme{\vPtl_\pm(T_\pm,\cdot)}^2_2 \nonumber.
    \end{align}
    For the rest of the proof, we fix $\gl = C_1^{-1}\max\ense{C_2 + C_1\lrpa{1+\norme{A}^2_\infty}; 2C_2}$, in order to have the following inequalities:
    $$  \begin{dcases*}
        C_1\gl^2 - C_2 \gl -C_2\lrpa{1+\norme{A}^2_\infty} &> 0, \\
        C_1\gl - C_2 &> 0,
    \end{dcases*}
    $$
    with the current value of $C_1$ and $C_2$.
    These last two inequalities and \eqref{ICa10}, allow us to find $\tau_1 > 0$ depending only on $\gO$, $T$, $\psi$ and $\norme{A}_\infty$ such that for every $\tau > \tau_1$, we have
    \begin{align}
        \label{ICa11}
        \norme{P^A_{\pm,\tau,\gl} \vPtl_\pm}^2_2 \geq& C_1 \lrpa{ \gl^2\tau^2 \norme{\vPtl_\pm}^2_2 + \gl\tau \norme{\nabla \vPtl_\pm}^2_2 + \tau \integdd{(\devp{\vPtl_\pm}{\nu})^2\abs{\devp{\psi}{\nu}}}{\gS_{\pm,\psi}}{\gs(x)}{t} } \nonumber \\
        &- C_2\tau \integdd{(\devp{\vPtl_\pm}{\nu})^2\abs{\devp{\psi}{\nu}}}{\gS_{\mp,\psi}}{\gs(x)}{t}.
    \end{align}
    \new

    Let us conclude this proof with the derivation of a suitable inequality similar to \eqref{ICa11} concerning $\vPt_\pm$ instead of $\vPtl_\pm$.
    We define $M = 2\max_\ade{\gO} \abs{\psi} > 0$ and notice that 
    $ e^{-\gl M} \leq e^{\pm 2\gl\psi} \leq e^{\gl M}$ on $\gO$.
    Recalling that $\vPtl_\pm = e^{\pm \gl\psi}\vPt_\pm $, we first get that $\nabla \vPtl_\pm = e^{\pm \gl\psi} \lrbr{\nabla \vPt_\pm \pm \gl \vPt_\pm \nabla\psi}$. It follows that 
    $$ \abs{\nabla \vPtl_\pm}^2 = e^{_\pm 2\gl\psi} \lrbr{ \abs{\nabla \vPt_\pm }^2 \pm 2\gl \vPt_\pm \nabla\psi \cdot \nabla \vPt_\pm + \gl^2 \vPt_\pm^2} \geq e^{-\gl M}\abs{\nabla \vPt_\pm}^2 \pm \gl e^{\pm 2\gl\psi}\nabla \psi \cdot \nabla \vPt_\pm^2, $$
    and by integrating by parts over $Q$,
    \begin{align*}
        \norme{\nabla \vPtl_\pm}_2^2 &\geq e^{-\gl M}\norme{\nabla \vPt_\pm}_2^2 \pm \gl \integdd{e^{\pm 2\gl\psi}\nabla \psi\cdot \nabla \vPt_\pm^2}{Q}{x}{t} \\
        &\geq e^{-\gl M}\norme{\nabla \vPt_\pm}_2^2 \mp \gl\integdd{\vPt_\pm^2 \dive\lrpa{e^{\pm 2\gl\psi} \nabla \psi}}{Q}{x}{t} \\
        &\geq e^{-\gl M}\norme{\nabla \vPt_\pm}_2^2 \mp \gl\integdd{\vPt_\pm^2 e^{\pm 2\gl\psi}\gD \psi}{Q}{x}{t} - \gl\integdd{\vPt_\pm^2 e^{\pm 2\gl\psi}2\gl\nabla\psi \cdot \nabla \psi}{Q}{x}{t} \\
        &\geq e^{-\gl M}\norme{\nabla \vPt_\pm}_2^2 - \gl\norme{\gD \psi}_\infty\integdd{\vPt_\pm^2 e^{\pm 2\gl\psi}}{Q}{x}{t} - 2\gl^2\integdd{\vPt_\pm^2 e^{\pm 2\gl\psi}}{Q}{x}{t} \\
        &\geq e^{-\gl M}\norme{\nabla \vPt_\pm}_2^2 - e^{\gl M}(\gl\norme{\gD \psi}_\infty + 2\gl^2)\norme{\vPt_\pm}_2^2.
    \end{align*}
    Moreover, we have
    $$ e^{-\gl M}\norme{\vPt_\pm}^2_2 \leq \norme{\vPtl_\pm}_2^2, \quad e^{-\gl M}\norme{\vPt_\pm(T_\pm,\cdot)}^2_2 \leq \norme{\vPtl_\pm(T_\pm,\cdot)}_2^2, $$
    along with, 
    $$ e^{-\gl M}\lrpa{\devp{\vPt_\pm}{\nu}}^2 \leq \lrpa{\devp{\vPtl_\pm}{\nu}}^2 \leq e^{\gl M}\lrpa{\devp{\vPt_\pm}{\nu}}^2 \quad \text{ and }\quad  \norme{P^A_{\pm,\tau,\gl} \vPtl_\pm}_2^2 \leq e^{\gl M} \norme{P_{\pm,\tau} \vPt_\pm }_2^2, $$
    since $P^A_{\pm,\tau,\gl} \vPtl_\pm = e^{\pm \gl\psi}P_{\pm,\tau} \vPt_\pm$.
    We then deduce from the lasts inequalities and \eqref{ICa11}, that for every $\tau > \tau_1$, 
    \begin{align*}
        C_1&\lrbr{ \lrpa{\tau^2-e^{2\gl M}(\norme{\gD \psi}_\infty + 2\gl)\tau}\gl^2 \norme{\vPt_\pm}^2_2 + \gl \tau \norme{\nabla \vPt_\pm}^2_2 + \tau \integdd{(\devp{\vPt_\pm}{\nu})^2\abs{\devp{\psi}{\nu}}}{\gS_{\pm,\psi}}{\gs(x)}{t} } \\
        &\leq e^{2\gl M}C_2\lrbr{ \norme{P_{\pm,\tau}\vPt_\pm}^2_2 + \tau \integdd{(\devp{\vPt_\pm}{\nu})^2\abs{\devp{\psi}{\nu}}}{\gS_{\mp,\psi}}{\gs(x)}{t} }.
    \end{align*}
    By taking $\tau_0 = \max\ense{\tau_1,2 e^{2\gl M}(\norme{\gD \psi}_\infty + 2\gl)}$, which depends only on $\gO$, $T$, $\psi$ and $\norme{A}_\infty$, we get for every $\tau > \tau_0$,
    \begin{align*}
        e^{2\gl M}C_2&\lrbr{ \norme{P_{\pm,\tau}\vPt_\pm}^2_2 + \tau \integdd{(\devp{\vPt_\pm}{\nu})^2\abs{\devp{\psi}{\nu}}}{\gS_{\mp,\psi}}{\gs(x)}{t} } \\
        &\geq C_1\lrbr{\gl^2\tau^2\norme{\vPt_\pm}^2_2 + \gl \tau \norme{\nabla \vPt_\pm}^2_2 + \tau \integdd{(\devp{\vPt_\pm}{\nu})^2\abs{\devp{\psi}{\nu}}}{\gS_{\pm,\psi}}{\gs(x)}{t} }.
    \end{align*}
    Finally, by taking $C = \frac{e^{2\gl M}C_2}{C_1 \min\ense{1,\gl,\gl^2}}$, which is a constant that depends only on $\gO$, $T$, $\psi$ and $\norme{A}_\infty$, we get \eqref{ICa2} which concludes this proof.
\end{Preuve}

Applying Theorem~\ref{Theo:Ca1}, we can also prove the following result:

\begin{Corollaire}
    \label{Cor:Ca3}
    For every $\vL_\pm \in H^{1,2}(Q)$ satisfying ${\vL}_{\pm|_\gS} = 0$, $\vL_{+|_{t=0}} = 0$ and $\vL_{-|_{t=T}} = 0$, there exist $\tau_0 > 0$ and $C >0$ which depend only on $\gO$, $T$, $\psi$, $\norme{q}_\infty$ and $\norme{A}_\infty$ such as for every $\tau > \tau_0$, we have
    \begin{multline}
        \tau \integdd{\abs{e^{\mp\gt_\tau}\devp{\vL_\pm}{\nu}}^2 \abs{\devp{\psi}{\nu}}}{\gS_{\pm,\psi}}{\gs(x)}{t} + \tau^2 \integdd{\abs{e^{\mp\gt_\tau}\vL_\pm}^2}{Q}{x}{t} + \integdd{\abs{e^{\mp\gt_\tau}\nabla \vL_\pm}^2}{Q}{x}{t} \\ \leq C\lrbr{ \integdd{\abs{e^{\mp\gt_\tau} L_\pm \vL_\pm}^2}{Q}{x}{t} + \tau \integdd{\abs{e^{\mp\gt_\tau}\devp{\vL_\pm}{\nu}}^2 \abs{\devp{\psi}{\nu}}}{\gS_{\mp,\psi}}{\gs(x)}{t} }.
    \end{multline}
\end{Corollaire}

\begin{Preuve}
    By taking $\vPt_\pm = e^{\mp\gt_\tau}\vL_\pm$, we get 
    $$ e^{\mp\gt_\tau}\nabla \vL_\pm = \nabla \vPt_\pm \pm \tau \vPt_\pm\nabla\psi, $$
    from which we get $\norme{e^{\mp\gt_\tau}\nabla \vL_\pm}_2^2 \leq 2 \lrbr{\norme{\nabla \vPt_\pm}^2_2 + \tau^2 \norme{\vPt_\pm}_2^2}$ and for every $\tau > 1$,
    \begin{multline*}
        \tau \integdd{\abs{e^{\mp\gt_\tau}\devp{\vL_\pm}{\nu}}^2 \abs{\devp{\psi}{\nu}}}{\gS_{\pm,\psi}}{\gs(x)}{t} + \tau^2 \integdd{\abs{e^{\mp\gt_\tau}\vL_\pm}^2}{Q}{x}{t} + \integdd{\abs{e^{\mp\gt_\tau}\nabla \vL_\pm}^2}{Q}{x}{t} \\ 
        \leq 2\lrbr{\tau^2 \norme{\vPt_\pm}^2_2 + \tau \norme{\nabla \vPt_\pm}_2^2 + \tau \integdd{\abs{\devp{\vPt_\pm}{\nu}}^2 \abs{\devp{\psi}{\nu}}}{\gS_{\pm,\psi}}{\gs(x)}{t}}.
    \end{multline*}
    We conclude this proof by using Theorem~\ref{Theo:Ca1}.
\end{Preuve}

\subsection{On the construction of the gauge function:}            \label{ConsGaug}

For the sake of completeness, we dedicate this last subsection to results used in the construction of the gauge function $\gvp$ for the proof of Theorem~\ref{Theo1}.

\begin{Proposition}
    \label{Pro:Ca4}
    Let $A$ be a compactly supported continuous vector field over $\RR^n$ such that $\dif{A}{}{} = 0$ (in the sense of distributions) over $\RR^n$ and $y \in \RR^n$. 
    Then, the map $\gvp : \RR^n \longrightarrow \RR$ defined by 
    $$ \forall x \in \RR^n ,\quad \gvp(x) = \Integd{A(y + s(x-y))\cdot (x-y)}{0}{1}{s} ,$$
    belongs to $\regu{1}{\RR^n}$ with $\nabla \gvp = A$ and $\gvp(y) = 0$.
\end{Proposition}

\begin{Preuve}
    We first prove that $\nabla \gvp = A$ in the sense of distributions before concluding over the desired regularity using mollifiers. 
    A proof of the case $y = 0$ can be found in \cite[Lemma 4.2.]{kianDeterminationNoncompactlySupported2019} for $n=3$; however, the proof can easily be rewritten for $n \geq 3$.
    For an arbitrary chosen $y \in \RR^n$, we set $\phi(x) = \gvp(y+x)$ and $B(x) = A(y+x)$ for every $x \in \RR^n$ such that 
    $$ \forall x \in \RR^n ,\quad \phi(x) = \Integd{B(sx)\cdot x}{0}{1}{s} .$$
    Since $B$ is a compactly supported continuous vector field over $\RR^n$ such that $\dif{B}{}{} = 0$ over $\RR^n$ then the case $y = 0$ gives us $\nabla \phi = B$ along with $\nabla \gvp = A$ over $\RR^n$, in the sense of distributions. 
    Now, let us prove that $\gvp \in \regu{1}{\RR^n}$ in order to conclude. 
    Let us fix $\rho \in \reguc{\infty}{\RR^n}$ such that $\supp(\rho) \subset B(0,1)$, $\rho$ positive and 
    $$ \integd{\rho(x)}{\RR^n}{x}  = 1.$$
    We then define $\rho_m  = m^n\rho(mx)$ and $\gvp_m = \gvp \ast \rho_m$ for every $m \in \NN$. 
    First, we know that the sequence $(\gvp_m)_{m\in \NN}$ converges pointwise to $\gvp$ on $\RR^n$.
    Secondly, we have $\gvp_m \in \reguc{\infty}{\RR^n}$ with, for every $i \in \iense{1}{n}$,
    $$ \devp{\gvp_m}{x_i} = \gvp \ast \devp{\rho_m}{x_i} .$$
    However, since $\nabla \gvp = A$ in the sense of distributions, then we obtain
    $$ \devp{\gvp_m}{x_i} = \gvp \ast \devp{\rho_m}{x_i} = A_i \ast \rho_m, $$
    and then $\lrpa{\devp{\gvp_m}{x_i}}_{m\in\NN}$ converges uniformly on $\RR^n$ with limit $A_i$ thus giving us $\gvp \in \regu{1}{\RR^n}$ with $\nabla \gvp = A$ over $\RR^n$.
\end{Preuve}

\begin{Proposition}
    \label{Pro:Ca5}
    Let $A \in \regub{0}{0,T; \regu{0}{\RR^n}}^n$ and $y \in \RR^n$. 
    Then $\gvp : (t,x) \in \inteoo{0}{T} \times \RR^n \longrightarrow \RR$ defined by
    $$ \forall (t,x) \in \inteoo{0}{T} \times \RR^n ,\quad \gvp(t,x) = \Integd{A(t,y + s(x-y))\cdot (x-y)}{0}{1}{s} ,$$
    belongs to $\regub{0}{0,T; \regu{0}{\RR^n}}$.
\end{Proposition}

\begin{Preuve}
    For any compact set $K'$ of $\RR^n$, let $\norme{\cdot}_{K',\infty}$ be defined as in the proof of Proposition~\ref{Pro:Mt1}.
    First, we show that $\gvp(t,\cdot) \in \regu{0}{\RR^n}$ for every $t \in \inteoo0T$.
    Let us fix $t \in \inteoo{0}{T}$ and set 
    $$ a : \RR^n \times \inteff{0}{1} \ni (x,s) \longmapsto A(t,y+s(x-y))\cdot (x-y) \in \RR.$$
    We have $a(\cdot,s)$ continuous for every $s \in \inteff{0}{1}$.
    Moreover, we have for $x \in \RR^n$, $M_x \geq 0$ such that 
    $$  \forall x' \in \ade{B(x,1)} ,\quad  \norme{a(x',\cdot)}_{\infty} \leq M_x, $$
    from which we get $\gvp(t,\cdot)$ is continuous at $x$ along with the first result.
    \new

    Now, let us prove that $\gvp \in \regu{0}{0,T; \regu{0}{\RR^n}}$.
    Let us fix $t \in \inteoo0T$ and consider a compact $K$ of $\RR^n$. 
    Since $A \in \regub{0}{0,T; \regu{0}{\RR^n}}^n$, for a fixed $\gve > 0$, there exists $\eta >0$ such that for every $t' \in \inteoo0T$,
    $$ \abs{t-t'} < \eta \implies \norme{A(t,\cdot) - A(t',\cdot)}_{K',\infty} \lrpa{ \sup_K \psi + 1} < \gve , $$
    where $K'$ is the compact subset of $\RR^n$ defined by
    $$ K' = \ense{ y + s(x-y) \given (x,s) \in K \times \inteff{0}{1} } .$$
    We deduce that for every $x \in K$ and $t' \in \inteoo0T$ such that $\abs{t-t'} < \eta$, 
    $$ \abs{ \gvp(t,x) - \gvp(t',x) } \leq \Integd{\abs{ A(t,y+s(x-y)) - A(t',y + s(x-y)) }}{0}{1}{s} \lrpa{ \sup_K \psi } <\gve, $$
    thus giving for every $t' \in \inteoo{0}{T}$,
    $$ \abs{t-t'}< \eta \implies \sup_{x\in K} \gvp(t,x) - \gvp(t',x) < \gve. $$
    Recalling that $\gve >0$ is arbitrary chosen, we get that $\gvp \in \regu{0}{\inteoo0T; \regu{0}{\RR^n}}$.
    \new

    We now conclude by proving that $\gvp \in \regub{0}{0,T; \regu{0}{\RR^n}}$.
    Since $A \in \regub{0}{0,T; \regu{0}{\RR^n}}^n$ then for any compact $K$ of $\RR^n$, we have $M_K$ such that for every $t \in \inteoo0T$,
    $$ \norme{A(t,\cdot)}_{K,\infty} \leq M_K . $$
    By fixing such $K$ and setting $ K' = \ense{ y + s(x-y) \given (x,s) \in K \times \inteff{0}{1} } \subset \RR^n$, we get that $K'$ is compact with 
    $$ \abs{\gvp(t,x)} \leq \Integd{\abs{ A(t,y+s(x-y)) \cdot (x-y)}}{0}{1}{s}  \leq M_{K'} \lrpa{ \sup_K \psi } $$
    for every $(t,x) \in \inteoo0T \times K$. 
    We conclude using the generality of $K$.
\end{Preuve}

\clearpage


\printbibliography

\end{document}